%
\font\ninerm=cmr9  	\font\eightrm=cmr8  	\font\sixrm=cmr6
\font\ninei=cmmi9	\font\eighti=cmmi8	\font\sixi=cmmi6
\font\ninesy=cmsy9	\font\eightsy=cmsy8	\font\sixsy=cmsy6
\font\ninebf=cmbx9	\font\eightbf=cmbx8	\font\sixbf=cmbx6
\font\ninett=cmtt9	\font\eighttt=cmtt8
\font\nineit=cmti9	\font\eightit=cmti8
\font\ninesl=cmsl9	\font\eightsl=cmsl8
%
\catcode`@=11
\newskip\ttglue
\def\tenpoint{\def\rm{\fam0\tenrm}%
  \textfont0=\tenrm \scriptfont0=\sevenrm \scriptscriptfont0=\fiverm
  \textfont1=\teni \scriptfont1=\seveni \scriptscriptfont1=\fivei
  \textfont2=\tensy \scriptfont2=\sevensy \scriptscriptfont2=\fivesy
  \textfont3=\tenex \scriptfont3=\tenex \scriptscriptfont3=\tenex
  \def\it{\fam\itfam\tenit}%
  \textfont\itfam=\tenit
  \def\sl{\fam\slfam\tensl}%
  \textfont\slfam=\tensl
  \def\bf{\fam\bffam\tenbf}%
  \textfont\bffam=\tenbf \scriptfont\bffam=\sevenbf
   \scriptscriptfont\bffam=\fivebf
  \def\tt{\fam\ttfam\tentt}%
  \textfont\ttfam=\tentt
  \tt \ttglue=.5em plus.25em minus.15em
  \normalbaselineskip=12pt
  \let\sc=\eightrm
  \let\big=\tenbig
  \setbox\strutbox=\hbox{\vrule height8.5pt depth3.5pt width\z@}%
  \normalbaselines\rm}
 
\def\ninepoint{\def\rm{\fam0\ninerm}%
  \textfont0=\ninerm \scriptfont0=\sixrm \scriptscriptfont0=\fiverm
  \textfont1=\ninei \scriptfont1=\sixi \scriptscriptfont1=\fivei
  \textfont2=\ninesy \scriptfont2=\sixsy \scriptscriptfont2=\fivesy
  \textfont3=\tenex \scriptfont3=\tenex \scriptscriptfont3=\tenex
  \def\it{\fam\itfam\nineit}%
  \textfont\itfam=\nineit
  \def\sl{\fam\slfam\ninesl}%
  \textfont\slfam=\ninesl
  \def\bf{\fam\bffam\ninebf}%
  \textfont\bffam=\ninebf \scriptfont\bffam=\sixbf
   \scriptscriptfont\bffam=\fivebf
  \def\tt{\fam\ttfam\ninett}%
  \textfont\ttfam=\ninett
  \tt \ttglue=.5em plus.25em minus.15em
  \normalbaselineskip=11pt
  \let\sc=\sevenrm
  \let\big=\ninebig
  \setbox\strutbox=\hbox{\vrule height8pt depth3pt width\z@}%
  \normalbaselines\rm}
 
\def\eightpoint{\def\rm{\fam0\eightrm}%
  \textfont0=\eightrm \scriptfont0=\sixrm \scriptscriptfont0=\fiverm
  \textfont1=\eighti \scriptfont1=\sixi \scriptscriptfont1=\fivei
  \textfont2=\eightsy \scriptfont2=\sixsy \scriptscriptfont2=\fivesy
  \textfont3=\tenex \scriptfont3=\tenex \scriptscriptfont3=\tenex
  \def\it{\fam\itfam\eightit}%
  \textfont\itfam=\eightit
  \def\sl{\fam\slfam\eightsl}%
  \textfont\slfam=\eightsl
  \def\bf{\fam\bffam\eightbf}%
  \textfont\bffam=\eightbf \scriptfont\bffam=\sixbf
   \scriptscriptfont\bffam=\fivebf
  \def\tt{\fam\ttfam\eighttt}%
  \textfont\ttfam=\eighttt
  \tt \ttglue=.5em plus.25em minus.15em
  \normalbaselineskip=9pt
  \let\sc=\sixrm
  \let\big=\eightbig
  \setbox\strutbox=\hbox{\vrule height7pt depth2pt width\z@}%
  \normalbaselines\rm}

%
%
\def\footnote#1{\edef\@sf{\spacefactor\the\spacefactor}#1\@sf
      \insert\footins\bgroup\eightpoint
      \interlinepenalty100 \let\par=\endgraf
        \leftskip=\z@skip \rightskip=\z@skip
        \splittopskip=10pt plus 1pt minus 1pt \floatingpenalty=20000
        \smallskip\textindent{#1}\bgroup\strut\aftergroup\@foot\let\next}
\skip\footins=12pt plus 2pt minus 4pt 
\dimen\footins=30pc 
\catcode`@=12  
%
%
\def\makearrow#1 #2 #3 #4 {%
   \setbox0=\hbox to #4pt{\rightarrowfill}%
   \wd0=0pt \ht0=0pt \dp0=0pt%
   \rlap{\hskip #1pt \raise #2pt%
      \vbox{
         \box0
         }%
      }%
   }
\input epsf
\font \trm = cmr12	
\font \bfsl = cmbxsl10  
\font \brm = cmr6	



\font \sym = msam10

	
\font \smc = cmcsc10	
\parindent=12 pt
\parskip = 0 pt

\def\donote#1#2{\footnote{${}^{#1}$}{#2}}
\def\section#1\par{\bigskip
  \message{#1}\centerline{\bf#1}\nobreak\smallskip\noindent}
\def\twolinesect#1#2\par{\bigskip
  \message{#1}\centerline{\bf#1}\nobreak
  \message{#2}\centerline{\bf#2}\nobreak\smallskip\noindent}
\def\threelinesection#1#2#3\par{\bigskip
  \message{#1}\centerline{\bf#1}\nobreak
 \message{#2}\centerline{\bf#2}\nobreak
\message{#3}\centerline{\bf#3}\nobreak\smallskip\noindent}
\def\subsection#1\par{\bigskip\noindent\message{#1}{\bf#1}\smallskip\noindent}

\def\therefore{.$\raise5pt\hbox{.}$.}
\def\sqr#1#2{{\vcenter{\vbox{\hrule height.#2pt
	\hbox{\vrule width.#2pt height#1pt \kern#1pt
		\vrule width.#2pt}
	\hrule height.#2pt}}}}


\def\lesnum{\hbox{\eightpoint{\ \raise2.5pt\hbox{$\prec$}\hskip-6.3pt\lower2.5pt\hbox{{\brm $\sim$}}\ }}}


\def\smoothaccent \raise3pt\hbox{{\sixrm '}}\hskip-4pt

\def\urltilde{\hskip1pt\lower3pt\hbox{$\tilde{\phantom { }}$ }}

\def\dia{{\sym\char'006} }
\def\diaa{{\sym\char'006}$\!$}
\def\diaas{{\sym\char'006}\nobreak\hskip -0.5 pt s }
\def\diab{{\sym\char'007} }
\def\diabb{{\sym\char'007}$\!$}
\def\diabs{{\sym\char'007}\hskip -0.5 pt s }

\hsize = 4.5 in
\vsize = 7.125 in
\hoffset = 1 in
\voffset = 1 in


\noindent\centerline{\trm Isonemal Prefabrics with Perpendicular Axes of Symmetry}
\vskip  5 pt
\centerline{R.S.D.~Thomas}

\centerline{St John's College and Department of Mathematics,}

\centerline{University of Manitoba, Winnipeg, Manitoba\ \ R3T 2N2\ \  Canada} 

\centerline{thomas@cc.umanitoba.ca\donote {*}{Work on this material has been done at home and at Wolfson College, Oxford.
Will Gibson made it possible to draw the diagrams with surprising ease from exclusively keyboard input. 
Richard Roth helped with the understanding of his papers and with this presentation.
An anonymous referee was most helpful.
To them all I make grateful acknowledgement.}
}
\vskip 5 pt 
\noindent {\bf Abstract} This paper continues the refinement of Richard Roth's taxonomy of isonemal weaving designs through types 11--32 of the 39 in order to solve three problems for those designs: which designs exist in various sizes, which prefabrics can be doubled and remain isonemal, and which can be halved and remain isonemal.
\vskip 5 pt
\noindent {\bf 1. Introduction} 

\noindent Except for a short list of interesting exceptions, Richard Roth [1993] has classified isonemal periodic prefabric designs into 39 infinite sets.
The first 10 have reflection or glide-reflection symmetries with parallel axes and no rotational symmetry. 
The next 22 have reflection or glide-reflection symmetries with perpendicular axes, hence half-turns, but no quarter-turns. 
The remaining 7 have quarter-turn symmetries but no mirror or glide-reflection symmetries. 
This paper is intended to reconsider this taxonomy for the designs with half-turns but no quarter-turns, refining it slightly to make it easier to use, and then to use it to examine three questions about prefabrics already discussed for what I call species 1--10 in a previous paper [Thomas, MS]. 
This reconsideration for species 11--32 presupposes the validity of Roth's taxonomy, which Roth has proved. 
The aim is to make it appear natural and to show it to be useful.
The answers to two of the questions---in sections 11 and 12---require this level of the taxonomy's refinement.

As Roth observes beginning his second paper [1995] on perfect colourings, `[r]ecent mathematical work on the theory of woven fabrics' begins with [Gr\"unbaum and Shephard, 1980], which remains the fundamental reference.
Roth's papers [1993] and [1995], however, contain the major advance from the fundamental work of Gr\"unbaum and Shephard.
In them he determines the various (layer, similar to crystallographic) symmetry-group types that periodic isonemal fabrics---actually prefabrics---can have and which of them can be perfectly coloured by striping warp and weft.
We are not concerned with striping, but the other terms and more are defined in [Thomas, MS], to which reference needs to be made.
\smallskip\noindent\epsffile{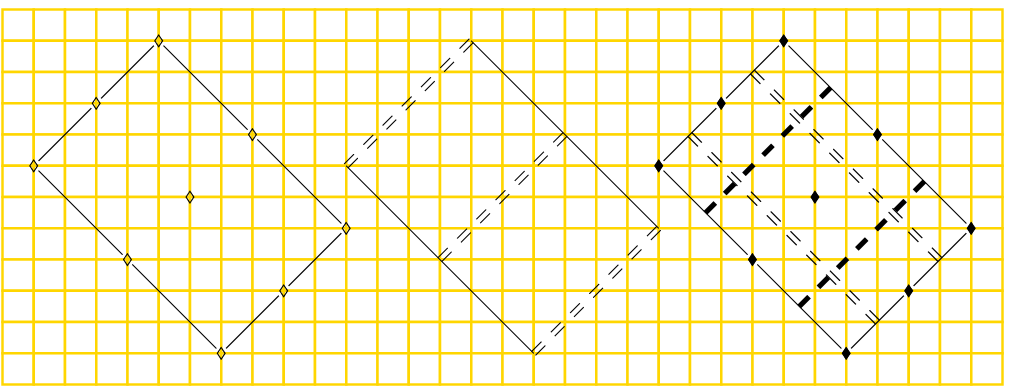}

\noindent\hskip 0.62 in (a)\hskip 1.16 in (b)\hskip 0.98 in (c)
\smallskip
\noindent\epsffile{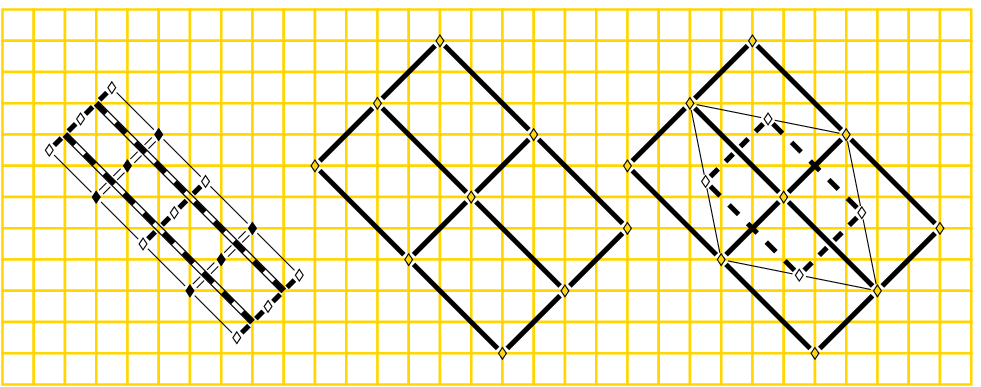}

\noindent\hskip 0.62 in (d)\hskip 1.16 in (e)\hskip 0.98 in (f)
\smallskip

\noindent Figure 1. The diagrammatic appearance of one lattice unit of some side-preserving subgroups $H_1$ and planar groups $G_1$, showing only short lengths of axes. \hskip 10 pt a. $H_1$ of type $p2$ with an outline of no geometrical significance through centres of half-turns. \hskip 10 pt b. $G_1$ or $H_1$ of type $pg$ with two boundaries that are axes of side-preserving glide-reflections and two conventional boundaries of no geometrical significance. \hskip 10 pt c. $G_1$ of type $pgg/pg$ with axes of side-preserving and side-reversing glide-reflections. \hskip 10 pt d. $H_1$ for $G_1$ of type $pmg/pgg$ with both sorts of glide-reflections and mirrors that are also axes of side-preserving glide-reflection. Two units of $G_1$ are shown.\hskip 10 pt e. $G_1$ of type $pmm/p2$. \hskip 10 pt f. $G_1$ of type $cmm/p2$ with four centres of half-turns and four mirrors outside the rhombic lattice unit.
\smallskip

{\it Lattice units}, that is, period parallelograms whose vertices are all images of a single point under the action of the translation subgroup.
Since our prefabric layers of strands meet at right angles and the symmetry groups with which we shall be concerned here all have parallel axes of reflection or glide-reflection in two perpendicular directions, the lattice units can be either rectangular or rhombic.
These rectangles and rhombs are defined by the group and representational conventions, as was not the case in [Thomas, MS].
The standard reference for symmetry terms is [Schattschneider, 1978] except that the term `order' is used here for one-dimensional period to distinguish it from two-dimensional period, which often differs.

For these prefabrics, the two-dimensional part $G_1$ of the symmetry group is of crystallographic type $pgg$, $pmg$, $pmm$, or $cmm$, that is, it has respectively glide-reflection axes in perpendicular directions, or glide-reflection axes and axes of reflection in directions perpendicular to each other, or just perpendicular axes of reflection, or both glide-reflection axes and axes of reflection (alternating) in both perpendicular directions.
In order for a glide-reflection to be a symmetry of the prefabric, it may or may not have to be combined with reversal of the sides of the prefabric, $\tau$. 
Mirror symmetry must always be reflection combined with $\tau$ because any cell through which the mirror passes must be its own image in the symmetry but the reflection alone, reversing warp and weft, would reverse its conventional colour. So $\tau$ is needed to restore it.
This means that there is never mirror symmetry in the side-preserving subgroup $H_1$. 
It is $H_1$ that determines the two-dimensional period under translation alone.
For these prefabrics, $H_1$ is of type either $p2$, which is a type of group generated by half-turns only (Figure 1a), or $pg$, generated by glide-reflections with parallel axes (Figure 1b), or $pgg$ explained in the first sentence of this paragraph.
When $H_1$ is of crystallographic type $pgg$ it may be the same group as $G_1$, or it may be a proper subgroup of $G_1$, in which case $G_1$ may not be of type $pgg$.
Possible geometrical configurations of the symmetry groups and side-preserving subgroups just mentioned and their lattice units, dimensions aside, are illustrated in Figures 1a to f, where the conventions for symmetry operations are displayed.
Glide-reflection axes without $\tau$ (side-preserving) are hollow dashed lines, while those with $\tau$ (side-reversing) are filled dashed lines.
Mirrors (always with $\tau$) are solid lines, normally filled, but with dashed filling when coincident with an axis of side-preserving glide-reflection.
Centres of half-turns are represented as \dia and of half-turns with $\tau$ as \diabb .
Thin lines are just boundaries of lattice units, outlining or completing the outline of rectangles or rhombs.
The lattice unit for type $pgg$ (Figure 1c) is displayed in combined side-preserving/side-reversing form to display the filled diamonds.
Finally, in Figure 1f four centres of half-turns at the intersection of mirrors are illustrated outside the type-$cmm$ $G_1$ lattice unit.
The longer and shorter diagonals of a rhomb will be called respectively its {\it length} and {\it width}.

As in species 1--10, none of the prefabrics discussed here has quarter-turn symmetry.
This fact has the consequence that each fabric mentioned here, with its specific orientation of lattice unit and axes within it has a quarter-turned image that is a distinct fabric, although not interestingly so.
This is true even for the initial subspecies $25_o$, $26_o$, $11_o$, $12_o$, $13_o$, and $14_o$ that have rectangular lattice units with both dimensions odd, since a $3\times 5$ rectangle is a $5\times 3$ rectangle rotated a quarter turn.

Since I am asking a reader to spend some time on a taxonomy, I ought to admit that no taxonomy that does not break its subject matter down into individuals is not subject to further refinement.
Richard Roth used something slightly coarser than his taxonomy to study perfect and chromatic colouring of fabrics because he did not need its refinement for that purpose.
In this paper I am determining a finer taxonomy because it is needed in \S\S 11 and 12, although not in \S 10. 
Different levels of refinement have different purposes; this one was designed primarily for \S 12.
Section 11 actually requires a small further {\it ad hoc} refinement of subspecies $17_e$, $19_e$, and $21$ that I make nothing of.

Successive sections will work through the species from 11 to 32 (\S\S 5, 6, 3, 7, and 9) with some organizational material first (\S 2), sections on rectangular and rhombic satins when they come up (\S\S 4 and 8), and at the end a section on each of the problems: doubling, halving, and which prefabrics exist with which orders (\S\S 10--12).

\smallskip
\noindent {\bf 2. Half-turns} 

\noindent The way I have chosen to discuss prefabrics with perpendicular axes of symmetry is to consider the half-turns first.
They and the type-$p2$ subgroups they generate will say a lot about the symmetry.
For one thing, these prefabrics are characterized by the presence of half-turns none of which are squares of quarter-turns.
(The absence of half-turns is the characteristic feature of prefabrics with only parallel symmetry axes (species 1--11), and quarter-turns are the characteristic feature of the remaining prefabrics, including species 33--39 and the exceptional prefabrics of order 4 and less.)
No prefabric in these species is of pure genus I or II [Gr\" unbaum and Shephard, 1985, p.~288; Roth, 1993, p.~202], and so all have half-turns to relate strands to adjacent strands, which relationship is characteristic of their isonemality.
(Half-turns can also appear not relating adjacent strands.)
The distribution of the half-turn centres relating adjacent strands is one of the main features of each symmetry group.
A lemma will be required.
\smallskip
\noindent {\smc Lemma 1.} The nine half-turn centres, with or without $\tau$, associated with the rectangular lattice unit of a symmetry group containing glide-reflections (types $pgg$ and $pmg$) lie either all at cell corners and perhaps cell centres or all on cell boundaries not at cell corners.
\smallskip
\noindent{\it Proof.} The fabrics illustrated in Figure 5 show that the former possibility is realized, and those in Figure 6 show that the latter possibility is realized.
What needs to be proved is that, if there are half-turn centres on cell boundaries not at a cell corner, then there are no half-turn centres at cell corners or cell centres.
The half-turn centres, \dia or \diabb , form a lattice in all $pgg$ and $pmg$ cases (Figures 1c and d).
The half-turn centres along each side of the rectangular lattice unit are taken by translation to those on the opposite side.
Accordingly, if one corner of the lattice unit is a half-turn centre on a strand boundary not at a corner, the other three are.
The sides of the lattice unit, being at $45^\circ$ to the strand directions, never pass through a cell corner or cell centre.
The four half-turn centres at mid-sides, which are halfway between these corners also lie on a strand boundary not at a corner.
The lines joining the opposite mid-side centres never pass through a cell corner or cell centre.
Finally, the half-turn centre at the centre of the lattice unit, which lies halfway between mid-side half-turn centres, lies on a strand boundary not at a corner.
The lemma is proved.
\smallskip
This paragraph is a description of how the half-turns relate to the strands and their symmetry; it should probably be skipped on first reading as it is heavily dependent on [Roth, 1993].
Some of the centres' distribution is determined by the induced strand subgroup of the symmetry group as follows:
each strand stands with respect to the half-turn centres in the same way by isonemality.
Using Roth's analysis (which should be consulted) and Roth's notation, analogous to $G_1/H_1$, for these induced strand subgroup types,

\noindent $\bullet$ $11/-$ indicates no symmetry of the strand except the translation by its order (species 11, 13, 15, 17, 19, 21, 25, 27, 29).
There must be \dia on one strand edge and either \dia or \diab (but not both) on the other edge.
The presence of these centres at cell corners has no implications for the distribution of dark and pale among the cells of the strand.

\noindent $\bullet$ If \dia and \diab appear on the same strand edge at cell corners, the composition of their half-turns effects a side-reversing translation necessarily by half the order, since its square is translation by the order.
This is the characteristic of $11/11$, which has both \diaas and \diabs along each strand edge (species 22, 30).

\noindent $\bullet$ Strands with reflective symmetry in the strand subgroup have centres of symmetry within the strand.
Note that \diab cannot occur inside a cell.
$12/-$ indicates no other symmetry  in the strand (species 12, 14, 18, 23, 26, 28, 31).
It happens in two ways.
Either there are \diaas in cells, in which case there must be also be \diaas or \diabs at cell corners along the edges so that there is strand-to-adjacent-strand transition (species $18_s$, 23, 26, $28_e$, $28_n$, 31), or there are \diaas between cells of the strand (species 12, 14, $18_o$, $18_e$), in which case there must also be \diaas on strand edges by isonemality.
It must be \dia on strand edges and not \diab because \diab on strand edges is characteristic of $12/11$ and $12/12$.
In one species there are \diaas at cell centres and corners and on edges too ($28_o$).
This is possible despite Lemma 1 because the symmetry is of type $cmm$.

\noindent $\bullet$ $12/11$ has side-reversing in-strand reflections that combine in pairs to give the order translation; \diab appears at cell boundaries, sides, top, and bottom (species 16, 20).

\noindent $\bullet$ $12/12$ has both \dia and \diab in both strands and their boundaries (species 24, 32).

However the half-turn centres are disposed with respect to each strand, a recurrent question is how they can be identically disposed with respect to every strand as required for isonemality.
To the answer to that question we now turn, using a pair of species of fundamental importance as examples.
\smallskip
\noindent{\bf 3. Crystallographic type} {\bfsl pmm}{\bf : Roth types 25 and 26} 

\noindent Every symmetry group of Roth types 11 to 32 has a specific subgroup of type $p2$ (Figure 1a), the nine half-turn centres of which are marked on a rectangle in Figures 1c to f (Figure 1f was expanded beyond a single lattice unit to include them).
Sometimes this subgroup is the side-preserving subgroup (Roth types 13, 14, 17, 18, 25, 26), but more often not.
(Fabrics of species 27 and 28 with symmetry of type $cmm$ have side-preserving subgroups of type $p2$ that are {\it not} those with the rectangular type-$p2$ lattice unit indicated but rather with the same rhomb as the $G_1$ lattice unit; type-$p2$ subgroups' lattice units do not need to be rectangular.)
How this configuration of nine centres is arranged will influence the action of each of these groups, and in turn how the other features of each group act will influence how the configuration is to be arranged.
\smallskip\noindent\epsffile{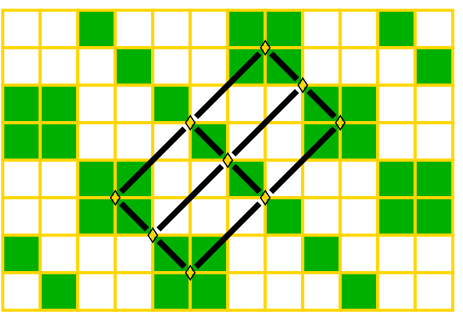}\hskip 10 pt \epsffile{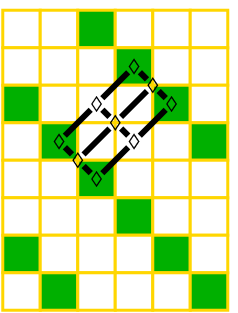}

\noindent (a)\hskip 1.8 in (b)
\smallskip

\noindent Figure 2. Roth's smaller examples of species 25 and 26. \hskip 10 pt a.~8-19-6%
\donote{1}{Such numbers, $a$-$b$-$c$, for fabrics refer to the catalogues of Gr\" unbaum and Shephard [1985; 1986]; a number $a$-$b$-$c^*$ refers to prefabrics that fall apart catalogued in [Hoskins and Thomas, 1991].} 
(species $25_e$). Four basic rectangles with their bounding mirrors and half-turns at the mirrors' points of intersection, always at cell corners. \break b. The 3/1 twill (species $26_e$). Four basic rectangles with their mirrors and half-turns in rows alternately of cell corners and cell centres.
\smallskip

\noindent {\it 25.} We begin with a minimum of constraints to see how the centres {\it might} be arranged by considering what may be the easiest type of symmetry group to understand, Roth type 25.
The crystallographic type of $G_1$ is $pmm$ illustrated in Figure 1e, and its side-preserving subgroup $H_1$ is what is left when the mirrors are removed, just half-turns (Figures 1a and 2a).
The pair is $pmm/p2$, which characterize species 25 and 26.
When the half-turns are restricted to cell corners, there are no implications for the strand symmetry, and so the strand subgroups are of type $11/-$, characteristic of type 25.
The question to which we now turn, and which will be the basis of consideration of all twenty-one other species, is how this configuration can be located among the cells of the design.
In particular, we need to see how the configuration can be located so as to make the fabric isonemal.
What is required is that there be half-turn centres on every one of two sets of parallel lines.
In this case the parallel lines are the strand boundaries,
and the strands are all related to their neighbours by them.
Since all the strands of the designs we are considering here must be so related (not the case in species 1 to 10), this condition is also necessary.
The $G_1$ lattice unit, because of the presence in it of two perpendicular mirrors, has built into it much redundancy; the four quarters are images of one another.
We need be concerned only about distribution of one quarter of a lattice unit to spread around the pattern, but also to spread around the half-turn centres. 
Having images of one quarter of some lattice unit placed so that a corner is on each of the parallel lines mentioned above is enough for isonemality.
The need to distribute some such rectangle, which we shall call a {\it basic rectangle}, is a feature of every one of these species.
As was the case for fabrics with only parallel symmetry axes, a basic rectangle will be so distributed to tessellate the plane when its dimensions, measured in diagonals $\delta$ of a cell, are relatively prime.
If the dimensions are $a\delta$ by $b\delta$---dimensions that will recur frequently (in Roth types 11, 13, 15, 17, 19, 21, 22, 23, 24, 27, 29, 30)---then the relative-primality condition $ca + db =1$ is a kind of recipe for how many basic rectangles need to be set out side by side in the two perpendicular directions to get a corner to the next line, vertical or horizontal.
To distribute the centres of half-turns to all strand boundaries, the basic rectangle (quarter of the lattice unit) must have this spacing, which we shall refer to as the {\it standard isonemal} spacing or constraint.
Roth's smaller example of species 25, 8-19-6 (Figure 2a), is smallest possible with $a=1$, $b=2$, and period 16.
When one parameter is odd and the other is even we shall call the species $25_e$, e for even.
When both are odd we shall call the species $25_o$, o for odd.
The while lattice unit has dimensions $2a\delta$ by $2b\delta$ and area $8ab$, the period of the fabric, whose order is $4ab$, half the period because the genus is III, and divisible by 4 in species $25_o$, divisible by 8 in species $25_e$.

\noindent {\it 26.} The complexity that is introduced by having half-turn centres not at cell corners will serve as an example of what will happen again and again as the simpler situation is extended to the less simple (roughly pairwise) throughout the remaining twenty Roth types.
If half-turn centres appear in cell centres (in all strands if in one) as well as at cell corners, then the strand subgroups are of type $12/-$, and the Roth type of the group is 26.
Along the mirrors, of which there are perpendicular parallel sets, there may be only in-cell centres or only cell-corner centres or alternating in-cell and corner centres (Figure 2b shows all three possibilities).
One set of parallel mirrors can contain alternately mirrors of in-cell and of corner centres, in which case the perpendicular mirrors contain exclusively in-cell and exclusively corner centres alternately (Figures 2b, 3a, and 4b). 
Or on the other hand the mirrors in both directions may have both in-cell and corner centres (Figures 3b and 4a).
\smallskip\noindent\epsffile{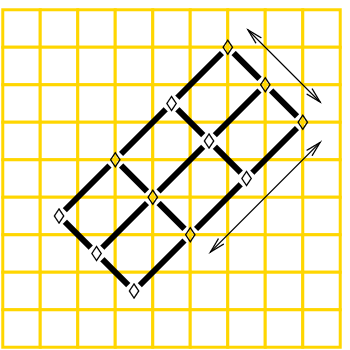}\hskip 10 pt \epsffile{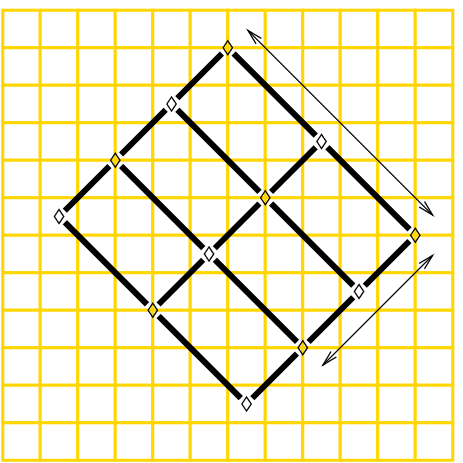}
\vskip -1.6 in\noindent\hskip 2.9 in $b\delta$
\vskip 0.1 in\noindent\hskip 1.15 in $b\delta$
\vskip 0.5 in\noindent\hskip 1.06 in $a\delta$
\hskip 1.8 in $a\delta$
\vskip 0.49 in\noindent (a)\hskip 1.34 in (b)
\smallskip


\noindent Figure 3. Arrangements of centres of half-turns at cell centres and cell corners in Roth type 26 showing the possibilities for alternative descriptions. The conventional description of these configurations is the upper 2/3 of each diagram rather than the lower 2/3.\hskip 10 pt a. Subtype $26_e$. \hskip 10 pt b. Subtype $26_o$.
\smallskip

\smallskip\noindent\epsffile{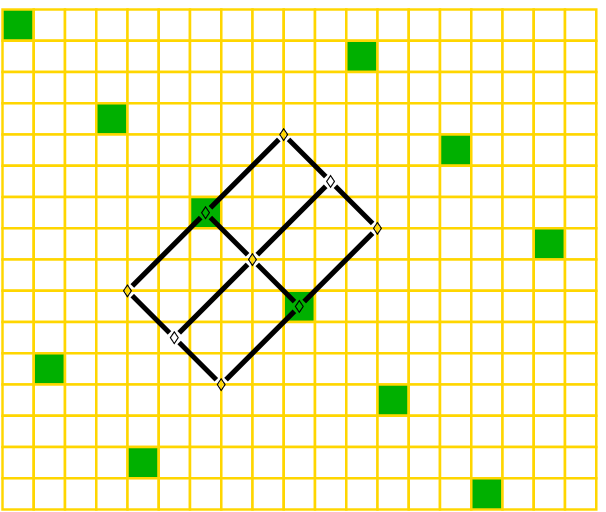}\hskip 10 pt \epsffile{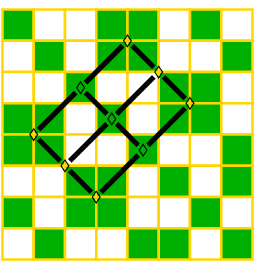}

\noindent (a)\hskip 2.35 in (b)
\smallskip

\noindent Figure 4. a. Species $26_o$ illustrated by the (30, 11) satin. One lattice unit with all mirrors carrying half-turn centres alternately at cell centres and cell corners. \hskip 10 pt b. Species $26_e$ illustrated by 12-619-1. One lattice unit with its half-turns in rows alternately of cell corners and cell centres.
\smallskip

These two possibilities are illustrated in Figure 3, and the illustrations show that the former case (Figure 3a) can be equally well described as having the configuration of the nine centres in a lattice unit be two lines of in-cell centres with a line of corner centres between them (lower 2/3) or be two lines of corner centres with a line of in-cell centres between them (upper 2/3). 
I shall use the second description.
Likewise the latter case (Figure 3b) can be described as having corners and centre of the lattice unit be in-cell centres (lower 2/3) or having corners and centre at cell corners, the other four centres being opposite (upper 2/3). 
I shall again use the second, so that in both diagrams I shall use the top 2/3 as the standard description.
Half the distance along an edge of one of these rectangles drawn on cells from cell corner to cell corner or cell centre to cell centre is $m\delta$ for some integer $m$,and between cell centre and cell corner is $\textstyle{(n+{1 \over 2})\delta}$ for some integer $n$.
The lattice-unit configuration with the cell centres down the middle (top 2/3 of Figure 3a) is characterized by being $a\delta$ by $b\delta$, where $b\delta$ is the length of the edge of cell corners (and line of cell centres) and so $a=2(n+{1 \over 2})$ is odd and $b=2m$ is even.
We shall call this subtype $26_e$.
The configuration with cell-corner corners and centre (top 2/3 of Figure 3b) has $a=2(n_1+{1 \over 2})$ and $b=2(n_2+{1 \over 2})$ both odd.
We shall call this subtype $26_o$.
It is easy to see why $a$ and $b$ cannot both be even; there would never be any cell-corner centres to serve as transitions between adjacent strands if lattice-unit corners were at cell centres (bottom 2/3 of Figures 3a and 3b).
Lattice-unit corners at cell corners with $a$ and $b$ even is characteristic of type 25.
Roth mentions [1993, p.~199] the 3/1 twill (Figure 2b) of species $26_e$ as an example of this type; it has $a=1$ and $b=2$, which is smallest possible. (Note that $3/1$ is not a fraction but the standard $2m$-tuple twill notation with $m=1$.)

To keep this Roth type from being the only one (among the 22 to be described here) with the dimensions of its whole $H_1$ lattice unit specified as $a\delta$ by $b\delta$, we introduce a smaller unit needed later $\beta \equiv \delta/2$, in which {\it one quarter} of the $H_1$ lattice unit, as discussed for species 25, is $a\beta$ by $b\beta$. 
These dimensions will occur again and again (Roth types 12, 14, 16, 20, 28, 31, 32).

Thinking along the lines of the discussion above of species 25, one could consider the distribution over the plane of one half of the subtype-$26_e$ lattice unit, the half with the even $b$ dimension divided.
One would then require that $(a,b/2)=1$; but this requirement does not differ from $(a,b)=1$ since $a$ is known to be odd.
So nothing is gained by dealing with $b/2$ instead of $b$, half the lattice unit instead of the whole.
We have no reason to care what power of 2 divides $b$.
In both species the period is $2ab$, as is the order, since the symmetry of the strands guarantees that a strand fragment of order length is a period parallelogram (rectangle).
For species $26_o$, $2ab$ is divisible by 2; for $26_e$, $2ab$ is divisible by 4.
The $(2ab-1)/1$ twill, among others, results from $a=1$, and the smallest non-twill examples must therefore have $a=3$, making the smallest order and period 12 in species $26_e$ and 30 in species $26_o$.
The former is confirmed by Roth's example 12-619-1 (Figure 4b), again of species $26_e$.
All of the rectangular satins, which have always been of interest, are of this type, and we pause to consider them.
\smallskip
\noindent{\bf 4. Rectangular satins} 

\noindent Rectangular satins are the simplest examples of species 26.
They are of even order [Gr\"unbaum and Shephard, 1980, Theorem 4] and a half-turn centre is in each dark cell, the others are in pale cells or at corners.
In species $26_o$, satins result from darkening cells containing half-turn centres along alternate mirrors.
In species $26_e$, satins result from darkening alternate cells containing half-turn centres along the mirrors containing just them in such a way that, along the perpendicular mirrors, such cells are either all dark or all not.
In Figure 3, these two verbal descriptions mean just to darken the corners of the alternative lattice units, the lattice units not marked with dimensions.
We have the condition in both subspecies $(a,b)=1$, equivalently $ca-db=\pm 1$ for some $c$ and $d$ (themselves with no common factor).
From the top corner of one of these alternatively described lattice units, $c$ lattice units down in the $a$ direction followed by $d$ lattice units back up in the $b$ direction will bring the corner of the unit to a strand one above or one below ($\pm 1$) its original level.
The horizontal distance moved, $ca+db$ cells, is the satin offset $s$ or $n-s$, where $n$ is the order.
Similarly $ba-ab=0$, and so $b$ units in the $a$ direction followed by $a$ units in the $b$ direction will bring the top (or any) corner to its original level $ba+ab$ units to the right of where it began.
The number of cells in the lattice unit $n = 2ab$, as we knew.

The algebraic condition of [Gr\"unbaum and Shephard, 1980, Theorem 4] for a rectangular satin is $s^2 -1 \equiv 0$ (mod $2n$), which is $(ca+db)^2 -1 \equiv 0$ (mod $4ab$) here.
Since $4abcd \equiv 0$ (mod $4ab$), and $ca-db = \pm 1$, 
$$0 \equiv (ca-db)^2 -1 + 4abcd = (ca+db)^2 -1 \hskip 24 pt(\hbox{mod}\ 4ab).$$
The smallest rectangular satin of species $26_e$ is the $(12,5)$ satin, but the smallest of species $26_o$ must be of order at least 30.
Since $11^2 -1 = 2 \times 60$, the $(30,11)$ satin is the smallest example (Figure 4a).
\smallskip
\noindent {\smc Theorem 1.} Rectangular satins are of two kinds, those of even order not divisible by 4 and those of order divisible by 4, falling in the species $26_o$ and $26_e$ respectively.
\smallskip
A satin results from distributing dark cells in the way characteristic of a satin, as in Figure 4a.
The simplest twill results from distributing one dark cell per order as a twill.
But other features of a prefabric can be distributed in such ways.
For instance, the cell-centre centres of half-turns in the fabric of Figure 4b have the distribution of the 5/1 twill and, less obviously, the cell-centre centres of half-turns in the fabric of Figure 4a have the distribution of the (15, 11) satin.
The (30, 11) satin of Figure 4a can be said to be based on the (15, 11) satin in the sense that the cell-centre centres of half-turns and the cell-corner centres of half-turns have congruent lattices with the (15, 11)-satin distribution.
Nothing is made to hang on this observation, but it can be a handy way to describe the distribution of half-turn centres.
It is just a consequence of isonemality that these underlying distributions must themselves be those of isonemal prefabrics.
In particular, it is not claimed that all features are distributed in such isonemal ways.
That such a claim is false is indicated by the example of the centres of side-preserving half-turns in genus-V prefabrics.
\smallskip
\noindent{\bf 5. Crystallographic type} {\bfsl pgg}{\bf : Roth types 11 to 16} 

\noindent We have considered the tessellation of the plane by lattice units for species 26 in terms of translation, but the actual way in which a single basic rectangle, which as in species 25 may not be a whole lattice unit, is distributed across the plane by the symmetry groups in question (of types $pgg$, $pmg$, and $cmm$) involve more than just translation and half-turns---glide-reflections and reflections to be precise, but the lattice units of all of those groups have associated with them rectangles that are distributed.
However they are distributed, the rectangles are those we have just considered for the rectangular satins and $(n-1)/1$ twills---the $(2m-1)/1$ twills to be precise ($n$ odd does not work) because those are the oblique rectangles that tessellate the plane in the regular way required.
I shall have to specify whether the inclusion of $a=1$ is allowed ($a$ rather than $b$ because $b$ may be required to be even but I adopt the convention that $a$ never is).

\noindent {\it 11.} The planar part of the symmetry groups of Roth types 11 to 16 is of crystallographic type $pgg$ with the glide-reflections side-preserving (the geometrical configuration of Figure 1c but with both sets of dashed lines hollow) or side-reversing (the geometrical configuration of Figure 1c with both pairs of parallel dashed lines filled) or of both kinds (Figure 1c).
All six types' situations can be analysed in terms of the simplest of them, that of Figure 5a (Roth type 11), where the glide-reflections are side-preserving so that the side-preserving subgroup is the whole group.
Roth's symbol is accordingly $pgg/-$.
With half-turn centres only at cell corners, there is no additional restriction on what sizes of basic rectangle with standard isonemal spacing can be the basis of this group's lattice unit.
The discussion above of moving the rectangle around was in terms of translations, and the movements here involve glide-reflection instead, but {\it where the rectangle goes} is the same in both.
The relevant basic rectangle is one quarter of the lattice unit, and the smallest possible basic rectangles with $a=1$, $b=2$, would, as lattice units, define the 3/1 twill.
They are the basis (in the sense of the final paragraph of \S 4) for Roth's smaller example of species $11_e$, 8-27-3.
His other type-11 example (with $a$ and $b$ both odd---subspecies $11_o$) is 12-111-2 based in this way on the 5/1 twill ($a=1$, $b=3$), but there is no need for the restriction to twills except to keep the order small.
Since the basic rectangle $a\delta$ by $b\delta$ is one quarter of the lattice unit, whose area is the two-dimensional period, the period is $8ab$ because the dimensions of the lattice unit are $2a\delta$ by $2b\delta$.
The first satin-based example of a group has $a=3$, $b=2$, $n=12$ and $s=5$.
This $n$ is the order of the satin that is the basis, but the order of the resulting fabric $4ab = 24$, half its two-dimensional period and always divisible by 4, by 8 if $b$ is even ($11_e$).
A fabric example with this group is illustrated in Figure 5a.
Since there are no implications for the internal symmetry of the strands, the strand symbol is $11/-$.
\smallskip\noindent\epsffile{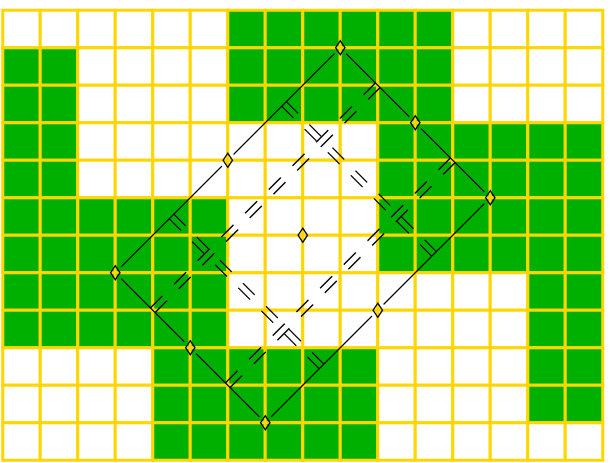}\hskip 5 pt \epsffile{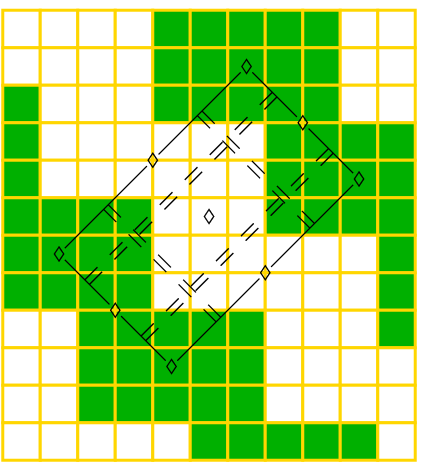}

\noindent (a)\hskip 2.3 in (b)
\smallskip


\noindent Figure 5. a. A smallest satin-based example of species $11_e$ based on the (12, 5) satin. \hskip 10 pt b. An example of species $12_o$ with basic rectangle 3 by 5 in $\beta$ units. Also an example of axes of glide-reflection not lying diagonally in cells, not in {\it mirror position,} in contrast to those in (a).
\smallskip

\noindent {\it 12.} When the half-turn centres can be off cell corners (Roth type 12), there are two possible sorts of location for the lattice unit.
The centres can be (a) a mixture of cell centres and cell corners or (b) all on cell boundaries (sides, tops, and bottoms with corners excluded).
Lemma 1 points out that the two classes of points cannot mix on lines at $45^\circ$ to the strand boundaries.
As a result, each of the groups that we consider will have to be examined to see whether it can be configured only the one way, only the other way, or both ways.
Hardly surprisingly, this simplest case can be configured in either position. 
Descriptions of both are eased by imagining a line down the middle of each strand.
The positions require half-turn centres to lie on these mid-lines, either (a) where the mid-lines cross one another (cell centres) and where strand boundaries cross (cell corners) or (b) where mid-lines cross strand boundaries, respectively.
We deal first with (a), then (b), then show that (a) and (b) exhaust the possibilities.

Consider placement of \dia at a cell corner.
If the next \dia in some direction is at a cell corner, then so will be all in that direction as successive images under the half-turns.
Since we are explicitly not considering the possibility that all are at cell corners, the next \dia in the perpendicular direction would have to be on a mid-line (i.e., at a cell centre); then alternating centres along that line (in both directions) will lie at cell centres.
And the same would be true of all parallel lines of centres, all being based at cell corners.
This configuration is impossible because neither a mirror nor a glide-reflection can relate a line of cell centres and a line of cell corners, and this configuration has alternating rows of each.
This configuration is impossible whenever axes of reflection or glide-reflection need to be placed {\it between} the lines of the nine-configuration (as in Figure 1d and 1f as well as the present 1c) for the same reason.

\smallskip\noindent\epsffile{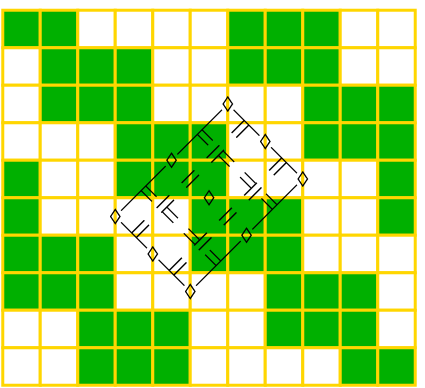}\hskip 10 pt \epsffile{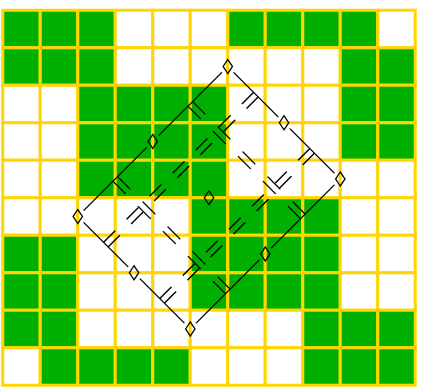}

\noindent (a)\hskip 1.64 in (b)
\smallskip

\noindent Figure 6. Examples of species $12_e$. a. 12-231-1 with $a=3$, $b=2$ in $\beta$ units.\hskip 10 pt b. Example with $a=3$, $b=4$ in $\beta$ units.
\smallskip

Placement of \dia at a cell corner must therefore have the next \dia in at least one direction at a cell centre, and so on alternately.
Consider that to be the case in one direction.
In the perpendicular direction, the positions can in principle either be maintained or reversed at every step.
To the first choice there is the objection of the previous paragraph, and so we find that if we have a cell corner at all we must have a checkerboard of cell corners and cell centres as in Figure 3b (ignoring its marked dimensions and mirrors) and Figure 5b.
In terms of $\beta$, the quarter lattice unit is an odd multiple $a\beta$ by an odd multiple $b\beta$.
For this reason, fabrics of this sort will be called species $12_o$ for odd.
Centres of both types will be distributed to all mid-lines and strand boundaries if and only if $(a,b)=1$, standard isonemal spacing.

This statement requires justification.
Since the basic rectangle has the half-turn centres in its adjacent corners in different positions in strands, and therefore in the same positions at diagonally opposite corners, perhaps the previous argument for $(a,b)=1$ does not apply.
Consider the four basic rectangles that compose the whole lattice unit $2a\beta$ by $2b\beta$.
The corners of the lattice unit are all in the same position in strands, and the geometrical aspect of the argument to and from $(a,b)=1$ works.
The arithmetical aspect delivers $2ca +2db = \pm 2$.
But this is what is wanted because we are counting strand widths by halves, and so a half-turn is moved to the corresponding point one strand from where it began, as required.

Neither $a$ nor $b$ can be 1 because that makes the design repeat so rapidly in the 1 direction that it is a twill, and all twills have mirror symmetry---unlike the designs we are considering.
We now turn to alternative (b) above.

For the sake of definiteness, consider placement of \dia at one side of a cell (as distinct from top and bottom).
All of the \diaas must therefore be on cell boundaries by Lemma 1.
Because a cell-side \dia is turned by both mirrors and glide-reflections into a \dia on the top or bottom of a cell, the images of the initial \diaa , diagonally opposite to it across basic rectangles making up the lattice units, are at top or bottom of cells while the \diaas at corners adjacent to it in those basic rectangles must likewise differ from each other for the same reason.
If the initial \dia is made the bottom corner of the lattice unit, the configuration of \diaas must be one of those in Figures 6a and 6b.
Once they are seen, they are seen to be just mirror images of each other (dimensions aside).
Next to the initial cell-side \dia in one direction or the other around any basic rectangle there is a cell-side \diaa , and parallel to that direction there are alternate rows of cell-side \diaas and top/bottom {\sym\char'006}s.
In terms of distance, the distance between two cell-side \diaas is an even multiple of $\beta$, and the distance between cell-side and top/bottom \diaas is an odd multiple of $\beta$.
For this reason fabrics of this sort will be called species $12_e$ for even.
Such spacings deliver \dia to every strand boundary and mid-line provided $(a,b)=1$.
They occur again in Roth types $14_e$, 16, $18_e$, and $20_e$.
Again $a=1$ would make the design a twill with a different symmetry, a symmetry with a group of this Roth type as a subgroup.

Consider placement of\  \dia at a cell centre.
The required strand-to-adja\-cent-strand transforms will be lacking unless there are also \diaas at cell boundaries or corners, and these possibilities have already been examined.
The division into (a) and (b) was indeed exhaustive.

The presence of \diaas on strand mid-lines gives the strand symbol $12/-$.
The smallest basic rectangles are $2\beta$ by $3\beta$ and $3\beta$ by $4\beta$ for species $12_e$ (side and top/bottom).
The first gives Roth's example 12-231-1 (Figure 6a); the second is illustrated in Figure 6b.
Order in general for species $12_o$ and $12_e$ is $2ab$.
The smallest basic rectangles are $3\beta$ by $5\beta$ for species $12_o$ (cell centres). 
As this gives the two-dimensional period and order 30, it is not surprising that such designs have not been noticed before.
An example is illustrated in Figure 5b.

\noindent {\it 13.} Two alterations to the two configurations so far considered need to be made to produce the remaining four Roth types with only glide-reflections.
The easier alteration is to introduce side-reversing (colour-preserving) glide-reflections in both directions.
The diagram is Figure 1c with both sets of parallel dashed lines filled (as in Figure 7a).
The type-11 discussion carries over completely in the geometrical part, but the side-reversing glide-reflections disappear from the side-preserving subgroup, making it be of type $p2$ but with the lattice unit of $G_1$.
The strand subgroups are of the same type, $11/-$, making the Roth type 13.
Roth's examples are 12-75-5 and 12-215-2 with basic rectangle $\delta$ by $3\delta$ (species $13_o$) because fabrics using $\delta$ by $2\delta$ (which would be in species $13_e$) have order 8 and have further symmetries.
Again there is no need for restriction to $\delta$ by $m\delta$.
An example of a fabric using $2\delta$ by $3\delta$ (order 24, period 48, species $13_e$) is illustrated in Figure 7a.
In general order is $4ab$, half the period, because the genus is III.
\smallskip\noindent\epsffile{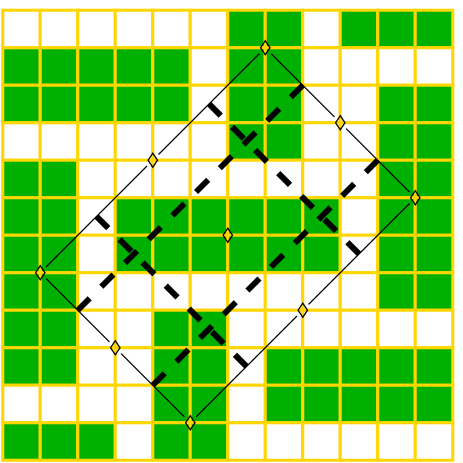}\hskip 10 pt \epsffile{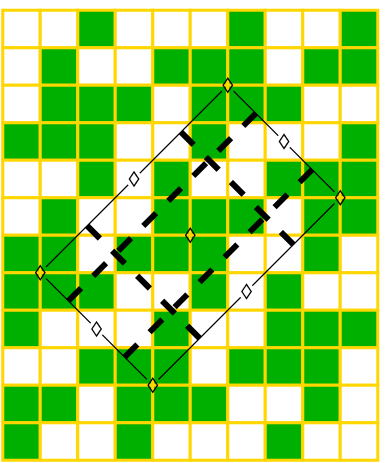}

\noindent (a)\hskip 1.8 in (b)
\smallskip


\noindent Figure 7. a.~An example of species $13_o$ based on the (12, 5) satin. \break
b. An example of species $14_o$ with basic rectangle 3 by 5 in $\beta$ units.
\smallskip

\noindent {\it 14.} The analogues of species $12_o$ and $12_e$ with side-reversing glide-reflec\-tions are species $14_o$ and $14_e$ with all argument the same as for species 12.
The outcome is different only in the side-preserving subgroup's type being $p2$ but with the $G_1$ lattice unit.
Examples again begin with smallest basic rectangles $2\beta$ by $3\beta$ for species $14_e$ (12-45-1) and $3\beta$ by $5\beta$ for species $14_o$, of which an example is illustrated in Figure 7b (order and period 30, a twillin with row-to-row offset 11---visible for warps in the figure). 
Order of both species $14_o$ and $14_e$ is $2ab$.

\noindent {\it 15.} The more interesting alteration to the type-11 scheme is to have side-preserving glide-reflection axes running one way and side-reversing glide-reflection axes running perpendicularly (Figure 1c).
This mixture does not change the two-dimensional geometry, but change there is.
The side-preserving subgroup has the one set of glide-reflections and not the other, making it be of type $pg$ (Figure 1b), but more interestingly the half-turns require $\tau$ for consistency.
Examples of what is species 15 work the same way as those of species 11, dividing into species $15_o$ and $15_e$, except that with $a=1$ and $b=2$, which is geometrically unproblematic, there is additional symmetry except for the prefabric that falls apart 8-5-2*.
Fabric examples accordingly begin with distances between \diab in $\delta$ units $a=1$ and $b=3$, which dimensions can be used both ways, as Roth shows with his examples: the long side of the basic rectangle (and lattice unit) can be parallel to the axes of side-preserving glide-reflection (12-47-5 in Figure 8a) or side-reversing glide-reflection (12-151-9 in Figure 8b).
The examples' order 12 is half of the period 24; in general order is $4ab$ (divisible by 4 for $15_o$ and by 8 for $15_e$).
Because of the \diabb , the fabrics are of pure genus IV.
\smallskip\noindent\epsffile{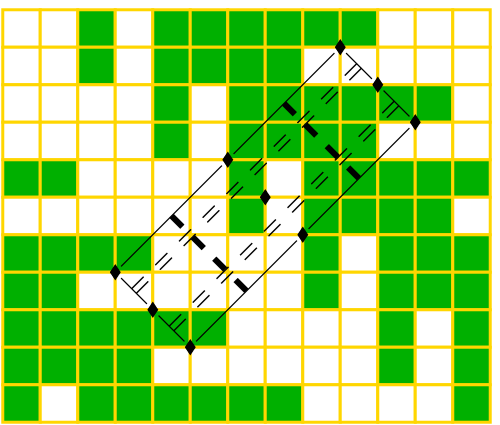}\hskip 10 pt \epsffile{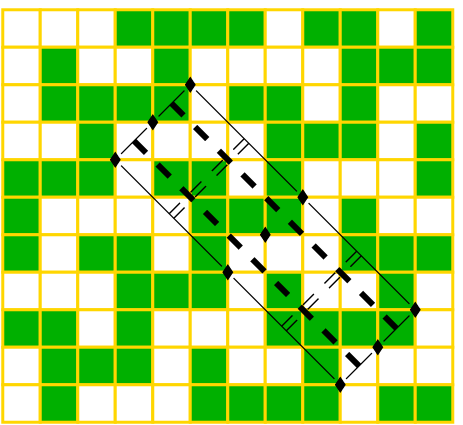}

\noindent (a)\hskip 1.92 in (b)
\smallskip


\noindent Figure 8. Roth's examples of species $15_o$. \hskip 10 pt a. 12-47-5. \hskip 10 pt b. 12-151-9.
\smallskip

The use of these rectangles `both ways' is an opportunity for a still finer distinction of subspecies into subsubspecies.
As remarked in the second-last paragraph of \S 1, sensible subdivision is relative to task.
Only once will we need to consider this particular sort of subdivision---for subspecies $17_e$ in the proof of Theorem 5.
I prefer with Roth [1993, p.~200] not to burden the reader elsewhere.

\noindent {\it 16.} When we turn to the symmetry groups with half-turn centres in the strands, we find that there is no analogue here, with the side-preserving and side-reversing glide-reflections, to species $12_o$ and $14_o$ because \diab cannot fall at cell centres.
Accordingly, with spacing analogous to that of $12_e$ and $14_e$ ($a>1$), species 16 has examples beginning with $a=3$, $b=2$, and period 12.
Order in general is $2ab$, divisible by 4.
Roth's examples (12-175-1, 12-603-1) are not only smallest possible but exhaust order 12 because other potential examples have additional symmetry.
The examples show that, as in species 15, lattice-unit dimensions can be used both ways.
The strand subgroups are different from those encountered before because the \diab in the strand gives half-turn symmetry only with side-reversal; without side-reversal, there is no symmetry except the basic translation.
Accordingly the strand subgroups are of type $12/11$.

A word on how period and order are related.
When a strand does not have half-turn symmetry within it, a fragment of it and either strand adjacent to it of order length make a period region (not necessarily a parallelogram), and period is twice order.
But when a strand has half-turn symmetry, then the second strand is not necessary (fabrics of genus I), and order length of one strand is a period parallelogram.
This accounts for the difference in relation between order and period, present but not commented on, between species 25, 11, 13, and 15 and species 26, 12, 14, and 16 respectively.
\smallskip
\noindent{\bf 6. Crystallographic type} {\bfsl pmg}{\bf : Roth types 17 to 24} 

\noindent The two-dimensional part of a symmetry group of types 17 to 24 is of crystallographic type $pmg$ with mirrors perpendicular to axes of glide-reflection---side-preserving or side-reversing or both.
Arbitrarily we begin with glide-reflections that are side-reversing; because mirrors are side-reversing too, the half-turns are side-preserving, as in Figure 9a and b.

\noindent {\it 17.} Whereas in Roth types 11--16 the axes of glide-reflection pass between the centres of half-turns at lattice-unit corners, in all type-$pmg$ symmetry the axes are boundaries passing through the the half-turn centres.
Since neither the side-reversing glide-reflections, which we consider first, nor the mirror reflections can be present in the side-preserving subgroup, only the half-turns remain (along with the translations), making it be of type $p2$.
If the half-turn centres are at cell corners only, so that their spacing is that of Roth type 11, there are no implications for symmetry of strands, the strand subgroups are of type $11/-$, and the Roth type is 17 (Figure 9a).
If, on the other hand, half-turn centres are allowed to fall at cell centres or mid-sides of cells, as in species 12, because of the action of mirrors on such centres (cell centre to cell centre, cell corner to cell corner, and cell top/bottom to cell side and vice versa) the centres fall either in alternating glide-reflection axes of cell centres only and of cell corners only (Figure 9b) or in glide-reflection axes containing alternating centres at cell top/bottom and at cell sides.
Either way yields strand subgroups of type $12/-$ and species 18.
\smallskip\noindent\epsffile{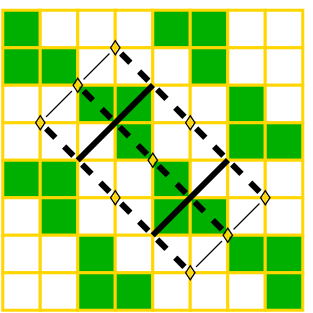}\hskip 5 pt \epsffile{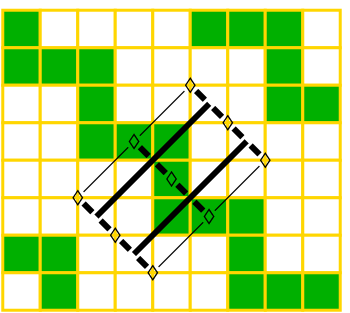} \hskip 5 pt \epsffile{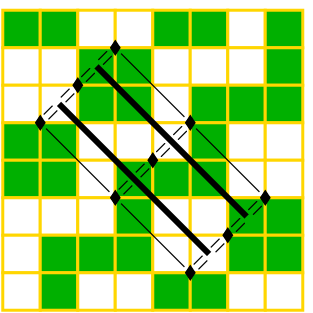}

\noindent (a)\hskip 1.11 in (b)\hskip 1.3 in (c)
\smallskip


\noindent Figure 9. a. 8-19-7 of species $17_e$. \hskip 10 pt b. 12-135-1 of species $18_s$. \hskip 10 pt c. 8-19-4 of species $19_e$.
\smallskip

Roth's examples of species 17 include 8-11-8 and 8-19-7 (Figure 9a), which illustrate that the same lattice unit ($2\delta$ by $4\delta$ with half-turn centres spaced like a 3/1 twill) can be used either way. 
The spacing can, however, be general rectangular-satin spacing, and so the quarter lattice unit can be odd multiple $a\delta$ by even multiple $b\delta$, starting with 1 by 2 as in Roth's examples of order 8 (species $17_e$), or odd by odd starting with 1 by 3 as in Roth's example of order 12, 12-11-2 (species $17_o$).
The condition $(a,b)=1$ is sufficient to produce isonemality in either species.
Order for both species is $4ab$, divisible as for species $11_o$ and $11_e$.

\noindent {\it 18.} Much of the spacing of half-turn centres in Roth type 18 is the same as in Roth type 12 (mirrors {\it were} mentioned there although there were none).
If half-turn centres lie alternately at cell top/bottoms and sides along axes, then the quarter lattice unit is an odd multiple $a\beta$ perpendicular to the mirrors by either an odd or even multiple $b\beta$ parallel to the mirrors since either will reach a line of tops/bottoms and mid-sides, giving two species $18_o$ and $18_e$ respectively.
If half-turn centres lie on alternate axes containing cell centres and cell corners, then one quarter of the lattice unit is an odd multiple $a\beta$ parallel to mirrors to change position (corner to centre) by an even multiple $b\beta$ perpendicular to mirrors to preserve position (centre to centre or corner to corner).
This quite different spacing has to be regarded as a third subspecies $18_s$.
Roth gives examples with the smallest possible lattice units of species $18_o$ and $18_s$, the 2/4 twill with $a=1$, $b=3$, and 12-135-1 (Figure 9b) with $a=3$, $b=2$, respectively.
If $a=1$ in species $18_s$, the fabric becomes a twill with mirrors perpendicular to the mirrors that are there; so $a=1$ is not allowed.
In species $18_e$, $a=1$ is allowed, however, but $b=2$ gives the 2/2 twill with too much symmetry.
But $a=1$ and $b=4$ give the 2/6 twill with the same symmetries as the 2/4 twill ($18_o$) but with a larger lattice unit and the minimum order 8 for species $18_e$.
Order is $2ab$ for all three subspecies, divisible by 4 in species $18_e$ and $18_s$, divisible by 2 in species $18_o$.

\noindent {\it 19.} The arbitrary choice that led to types 17 and 18 can now be reversed to combine side-preserving glide-reflections with mirror symmetry.
Because the one is side-preserving and the other side-reversing, the half-turns must be side-reversing as in types 15 and 16.
That in turn means that they disappear from the side-preserving subgroup along with the mirrors, making $H_1$ be of type $pg$.
If the half-turn centres are at cell corners only, so that the spacing is that of type 17, there are no implications for symmetry of strands, the strand subgroups are of type $11/-$, and the Roth type is 19.
If, on the other hand, half-turn centres are allowed to fall off cell corners, as in type 18, the fact that \diab cannot fall in a cell prevents the placement analogous to species $18_s$.
As in the discussion of type 18, the centres fall on glide-reflection axes containing alternating centres at cell top/bottom and at cell sides, the spacing being those of species $18_o$ and $18_e$, the strand subgroups being of type$12/11$ on account of in-strand \diabb, and Roth type being $20$.

In species 19, the smallest basic rectangle, measured in $\delta$ units, has $a=1$, $b=2$ and gives Roth's smaller example of species $19_e$, 8-19-4 (Figure 9c), with period 16.
Roth's larger example is of species $19_o$, 12-23-4 with $a=1$, $b=3$, and period 24.
General standard isonemal spacing is possible, and the lattice units can be used both ways (8-11-4 and 12-55-5 respectively use these units the other way).
The genus is pure IV, and the order is $4ab$ for both species, divisible by 4 in $19_o$ and by 8 in $19_e$.
As in species 15, the genus prevents $H_1$ from being transitive on strands.

\noindent {\it 20.} Roth illustrates both subspecies of species 20 with the smallest possible examples, the 2/2/1/1 twill having $a=1$, $b=3$, in $\beta$ units for the basic rectangle in species $20_o$, and 12-371-1 having $a=3$, $b=2$, in species $20_e$.
As in species $18_o$, in species $20_o$ lattice units can be used both ways once 1 is not a parameter.
Order for both species is $2ab$, with divisibility as for $18_o$ and $18_e$.

\noindent {\it 21.} The twist introduced in these Roth types is that the glide-reflections perpendicular to the mirrors need not be of the same kind.
It is possible to have alternating side-preserving and side-reversing glide-reflections as in species 3 and 4 (See Figures 10a and 11).
Along the axes of side-preserving glide-reflections the half-turns are side-reversing, while along the alternate (side-reversing) axes the half-turns are side-preserving.
The latter axes disappear from the side-preserving subgroup $H_1$ leaving behind the centres of half-turns along them while the former axes remain but lose their side-reversing {\sym\char'007}s, and when the (side-reversing) mirrors go they leave behind side-preserving glide-reflections so that $H_1$ is of type $pgg$.
The axes of side-reversing glide-reflection become the boundaries (and centre line) of the lattice unit for $H_1$, the dimension of which parallel to the original mirror direction is four times instead of two times the distance between the original axes; it can be regarded as two of the $G_1$ lattice units side by side\break
\smallskip\noindent\epsffile{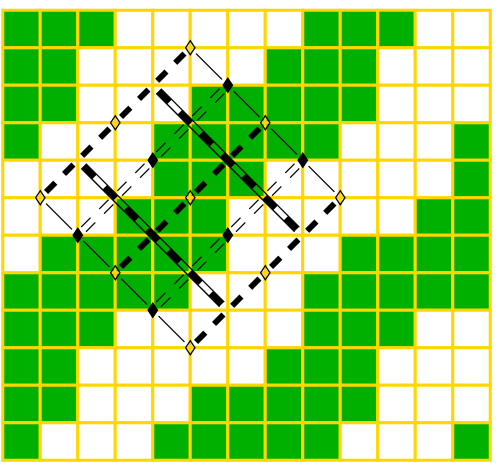}\hskip 10 pt \epsffile{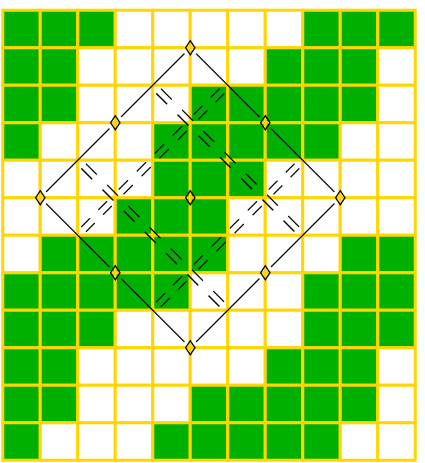}

\noindent (a)\hskip 1.92 in (b)
\smallskip


\noindent Figure 10. 8-7-2 of species 21. \hskip 10 pt a.~Two $G_1$ lattice units. \hfil\break
b.~The corre\-sponding $H_1$ lattice unit.
\smallskip

\noindent (See Figure 10b).
Thus it is reasonable to take an axis of side-preserving glide-reflection as the centre line of the $G_1$ lattice unit (a side-reversing axis {\it could} be used).
The rationale for this choice is the same as for the different-looking choice for species 3 and 4 in [Thomas, MS].
The illustrations in Figure 10 of 8-7-2 show two lattice units (so described) of $G_1$ and one of $H_1$.
If centres of half-turns occur only at cell corners and the \diaas and \diabs occur on alternate strand edges (never on the same edge), there are no implications for the symmetry of the strands so that the strand subgroups are of type $11/-$.
The distance parallel to mirrors between closest \diaas is automatically an even multiple of $\delta$ on account of the cell-corner \diab between them, and vice versa.
Thinking of parallel strand boundaries' being labelled with integers both ways from zero allows us to consider the even-numbered and the odd-numbered boundaries and to insist that the \diaas fall on the one set and the \diabs on the other.
This condition forces the distance perpendicular to mirrors between closest \diaas (and between closest {\sym\char'007}s) to be an even multiple of $\delta$.
It therefore also forces the distance parallel to mirrors between \dia and closest \diab to be odd (or all centres would fall on boundaries of the same parity).
The quarter $G_1$ lattice unit must therefore have standard isonemal spacing but with $a\delta$ having $a$ odd parallel to mirrors and $b\delta$ with $b$ even perpendicular to mirrors (no analogue of $17_o$ and $19_o$).
As the example 8-7-2 of Figure 10 shows, $a=1$ is allowed (there $b=2$).
This is the characteristic pair $pmg/pgg$, $11/-$, of Roth type 21, which is our first example of a type where the period is not the area of the lattice unit of $G_1$.
The period is determined by $H_1$, and here its lattice unit has twice the area of $G_1$'s.
In the example of Figure 10, the period is $32=16ab$ generally, four times the order, because a piece of strand of order length appears in four guises in the period region formed by transforming it three times to successively adjacent strands by half-turns alternately \dia and \diabb .
The order is $4ab$, divisible by 8.
The genus is pure V, and so $H_1$ is not transitive on strands.
The next three Roth types also deliver anomalous period calculations.

\noindent {\it 22.} If the centres of half-turns are restricted to cell corners but \dia and \diab are allowed to appear on the same edge of a strand---and so both on every strand edge, the strand subgroups are of type $11/11$, characteristic (with $pmg/pgg$) of Roth type 22.
This is guaranteed by having \dia fall on every strand boundary, achieved by giving standard isonemal spacing to one quarter of the $H_1$ lattice unit, which is one half (divided parallel to the mirrors) of the $G_1$ lattice unit.
But the distance parallel to mirrors between closest \diaas must still be even to fit a \diab between them, and so the distance between closest \diaas perpendicular to mirrors must be odd, making the quarter $H_1$ lattice unit be odd $a\delta$ (perpendicular to mirrors) by even $b\delta$ (parallel to mirrors), $(a,b)=1$.
That is not unusual, but the $G_1$ lattice unit is $2a\delta$ by $b\delta$.
Despite the even-by-even dimensions of the $H_1$ lattice unit, $H_1$ is transitive on the strands because \dia is distributed to every strand boundary.
Roth's examples of this species, 8-15-1 and 8-45-1, both have $a=1$, $b=2$, which is smallest possible.
In general order is $4ab$, divisible by 8.
The parameters need not be so small, but because the period is $8ab$, twice the order, it does grow quickly. 
With $a=3$, $b=2$, the next possibility, the spacing of all the half-turn centres is still that of a twill.
The first at non-twill spacing of all centres has $a=3$, $b=4$, order 48, and period 96; an example is illustrated in Figure 11a.
\smallskip\noindent\epsffile{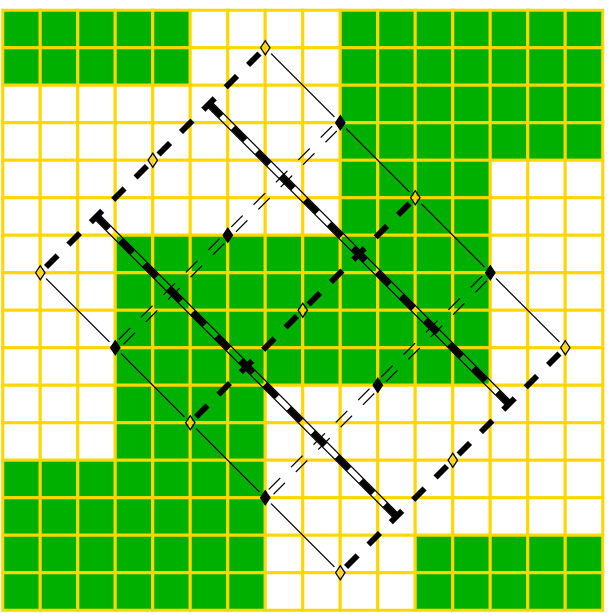}\hskip 5 pt \epsffile{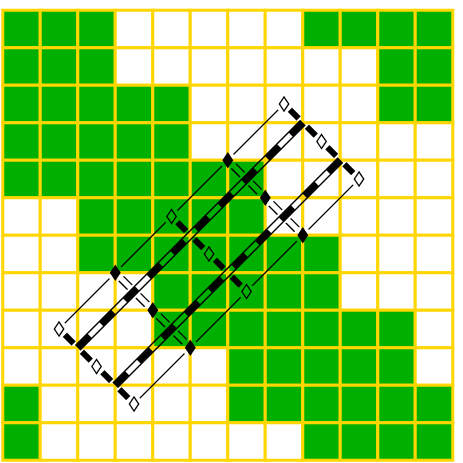}

\noindent (a)\hskip 2.3 in (b)
\smallskip


\noindent Figure 11. Pairs of $G_1$ lattice units composing $H_1$ lattice units. a. Example of species 22 of order 48. \hskip 10 pt b. 12-31-1 as an example of species $23_o$ ($a=3$, $b=1$).
\smallskip

\noindent {\it 23.} The simpler of the two ways to allow centres off cell corners is to allow \dia in cell centres and keep the \diabs at cell corners along axes in mirror position.
This will make the strand subgroups of type $12/-$ and the Roth type 23.
The rows of alternating \diaas and \diabs are still parallel to the mirrors (Figure 11b).
Because a cell-centre \dia falls between cell-corner {\sym\char'007}s, the distance between closest \diaas (or {\sym\char'007}s) parallel to the mirrors is an odd multiple $a\delta$.
As Roth's examples show, $a=1$ is allowed in species 23.
Distribution of \dia to all strands (\diab to all strand boundaries) is accomplished by having each half of the $G_1$ lattice unit (divided parallel to the mirrors---one quarter of the $H_1$ lattice unit as in species 22) have standard isonemal spacing.
As the dimension parallel to the mirrors is already assigned $a\delta$ with $a$ odd, the dimension perpendicular to the mirrors can be odd or even $b\delta$ so long as $(a,b)=1$.
Roth's examples of this species are of both kinds. 
While $a=1$ in both, with $b=2$ he has 8-27-1, and with $b=3$ he has 12-107-2 with periods 16 and 24 respectively.
For $a=1$, $b=1$, there is the prefabric that falls apart 4-1-1*, but $b$ can be 1 for larger $a$, as in 12-31-1 ($a=3$, period 24, Figure 11b), and both can be larger.
This division of $b$ into odd or even gives species $23_o$ and $23_e$.
For both species, order is $4ab$, divisible by 8 for $23_e$, by 4 for $23_o$.
As in species 22, the lattice unit of $H_1$ is even by even so that $H_1$ is not transitive on strands by translations, and the \diaas fall in cells not on strand boundaries; as a result $H_1$ is not transitive.

\noindent {\it 24.} The more complicated way to allow centres off corners is to have them all on cell tops/bottoms and sides as in species $12_e$, $14_e$, 16, $18_e$, and $20_e$.
This puts both \dia and \diab in each strand, giving strand subgroups of type $12/12$ and characterizing Roth type 24.
This configuration can be extremely compressed as in the smaller of Roth's examples, the 2/2 twill, where the order and period are only 4.
The dimension of the $G_1$ lattice unit across the mirrors must be odd $a\delta$ to accommodate the mirrors' interchanges of types of half-turn centre positions.
The \dia at the mid-point of the lattice-unit side perpendicular to the mirrors will be in the opposite position to that of the adjacent corners (top/bottom rather than side or vice versa).
The four corners of the $G_1$ lattice unit, however, will be in the same position in order to accommodate the \diabs on the sides parallel to the mirrors.
The \diab at the mid-point of those sides will be in the same position as or opposite position to that of the corners depending on whether the side length is odd or even.
The \diab at the centre of the lattice unit will be of the opposite type.
And the lattice unit has alternative descriptions since its boundaries can be moved perpendicular to the mirrors the distance between them.
There are four lattices, that of \dia on warp boundaries, that of \dia on weft boundaries, that of \diab on warp boundaries, and that of \diab on weft boundaries.
All behave similarly, distributing their centres to all parallel strand boundaries if any does.
One way to express the requirement is to have the \diaas at the corners of the $G_1$ lattice unit distributed to all parallel strand edges.
For that, the whole $G_1$ lattice unit must have standard isonemal (including $(2m-1)/1$-twill) spacing in $\delta$.
Since the dimension of the lattice unit across the mirrors is odd $a\delta$, the dimension of the $G_1$ lattice unit parallel to the mirrors can be odd or even $b\delta$, so long as $(a,b)=1$, as the examples of the 4/4 twill ($a=1$, $b=2$) and 6/6 twill ($a=1$, $b=3$) show.
Another way of expressing this condition is that the $G_1$ half of the $H_1$ lattice unit is $a\delta$ by $b\delta$, making the period and order $4ab$, divisible by 4 for $24_o$ ($b$ odd), by 8 for $24_e$ ($b$ even).
Odd pairs can be used either way.
Roth's larger example, 12-189-1 ($a=3$, $b=1$), is an example of a pair already used one way (in the 6/6 twill) being used in the other way.
\smallskip
\noindent{\bf 7. Crystallographic type} {\bfsl cmm} {\bf with side-reversing glide-reflec\-tions: Roth types 27 and 28}

\noindent {\it 27.} When additional symmetry is introduced to the configuration of Roth types 25 and 26 ($pmm$) but without creating the much greater (rotational) symmetry of a few exceptional prefabrics (crystallographic type $p4m$), the type is $cmm$.
The standard lattice unit of $cmm$ is only part of the configuration illustrated in Figure 1f.
The nine half-turn centres are still there, but the lattice unit is the rhomb with edges joining the half-turn centres in the middle of the nine-configuration's sides.
At the mid-points of the edges of this rhomb are {\it supplementary} half-turn centres at the intersections of glide-reflection axes.
The supplementary half-turn centres and such axes between the mirrors are what distinguish groups of type $cmm$ from those of type $pmm$.
Since the rhombs tessellate the plane, their corners must all have the same position in a cell.
For feasibility the length and width of the rhomb are of the same parity (Section 3 of [Thomas, MS]).
The presence of the rhomb places no other constraint on the placement of the configuration of the nine half-turn centres, but we begin with the simplest situation where all thirteen half-turn centres are at cell corners (a restriction that will gradually be lifted).
We also begin with the glide-reflections side-reversing and the half-turns that they generate side-preserving.
The side-preserving subgroup $H_1$ will accordingly contain all of the half-turns but none of the axes.
Its lattice unit is the rhomb.
The groups are of types $cmm/p2$, and because the half-turn centres are at cell corners there are no implications for strand symmetry, making the strand subgroups be of type $11/-$.
This combination is characteristic of Roth type 27.
\smallskip\noindent\epsffile{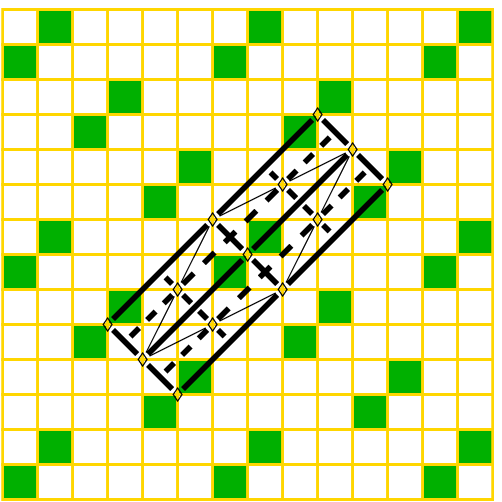}\hskip 5 pt \epsffile{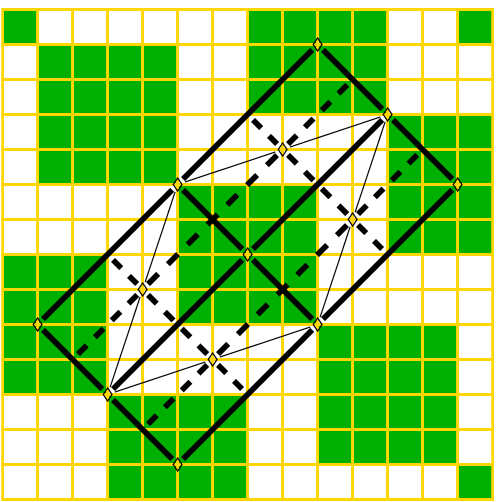}

\noindent (a)\hskip 1.86 in (b)
\smallskip


\noindent Figure 12. a. 6-1-1 as example of species $27_o$. \hskip 10 pt b. Order-16 example of species $27_e$.
\smallskip

We must consider what dimensions will make fabrics isonemal.
In the tessellation of the plane by the rhombs, their central rectangles, bounded by glide-reflection axes and with corners at the supplementary half-turn centres, cover the plane with gaps congruent to themselves.
Since it is their corners that need to be distributed, the gaps do not matter.
So it would distribute half-turn centres to all strand boundaries to have those central rectangles obey the isonemal constraint with odd-odd sides in $\delta$ ($a=1$ allowed).
Odd-even will not do here because it puts the half-turn centre at the centre of the rhomb not at a cell corner; this omission will be filled in.
The spacing just proposed is not necessary, though it is sufficient, because all of the work of strand-to-adjacent-strand transformation is being done by the supplementary half-turns and none by the usual nine.
To use the nine half-turns at the mirror intersections, the two lattices of half-turn centres that can be seen in Figure 12---those at mirror intersections and those at glide-reflection-axis intersections---must be offset not just by an odd multiple of $\delta$ as in Roth type 21, where lattices of \diaas and of \diabs had to be distributed both to strand edges but not both on the same side of a strand, but by having each \dia of one lattice be the centre of a rectangle of the other lattice so that each lattice's \diaas affect only alternate strand boundaries.
As in Roth type 21, the two lattices have the same spacing, the inner rectangle being congruent to a quarter of the nine-configuration.
The dimensions of the central rectangle (and the basic rectangle) may accordingly have greatest common factor twice those of the standard isonemal constraint.
Any common factor larger than two will cause strand boundaries to be skipped.
So these rectangles can have dimensions $a\delta$ and $b\delta$ or $2a\delta$ and $2b\delta$ with $(a,b)=1$.
If $a$ and $b$ are both odd in the $2a\delta$-by-$2b\delta$ case, then the two lattices fall on the same boundaries; so $a$ and $b$ must be odd and even.
These species can be called $27_o$ and $27_e$ respectively.
The orders are respectively $2ab$, divisible by 2, and $8ab$, divisible by 16.
The smallest are $\delta$ by $3\delta$ and $2\delta$ by $4\delta$ respectively.
Roth's examples are both of the former sort, 6-1-1 (Figure 12a) with $a=1$, $b=3$, period 12, and 10-19-2 with $a=1$, $b=5$, period 20.
An example of the latter sort with $2\delta$ by $4\delta$, $a=1$, $b=2$, period 32, and order 16 is illustrated in Figure 12b.%
\smallskip
\noindent {\it 28.} Roth type 28 is characterized by the same $G_1/H_1$ pair as 27, but the strand subgroups of type $12/-$ arise from \diaa s' falling at cell centres.
For the rhomb to exist the parities of its length and width, the dimensions of the nine-configuration, must be equal.
The mid-points of the rhomb sides need not fall at cell corners but can fall either there or at cell centres or on cell boundaries.
Since the four corners of the central rectangle are images in mirrors of one corner, they are all in the same sort of position, making three different ways of arranging these rectangles.
To produce type 27 it was necessary to keep half-turn centres away from cell centres and boundaries to achieve $11/-$; here it is possible to use cell centres to achieve $12/-$.

If the central rectangle has corners at cell corners and has standard isonemal spacing to distribute them, odd-odd dimensions in $\delta$ units make its centre and the other eight centres of the nine-configuration fall also at cell corners as in the previous paragraph. 
This must be ruled out.
Standard isonemal spacing and odd-even dimensions for the central rectangle, on the other hand---the central rectangle still having corners at cell corners---make the nine (other) half-turn centres all be cell centres and be appropriately distributed because a quarter of the nine-configuration has the same dimensions as the central rectangle.
Roth's explicit example of this type, the $(8, 3)$ satin 8-1-1 (Figure 13a), is of this kind with $a=1$, $b=2$.
We shall call this species $28_e$ for $b$ even 
and because it is exemplified by even-order rhombic satins.
One sees that its period and order are $4ab$ divisible by 8.
\smallskip\noindent\epsffile{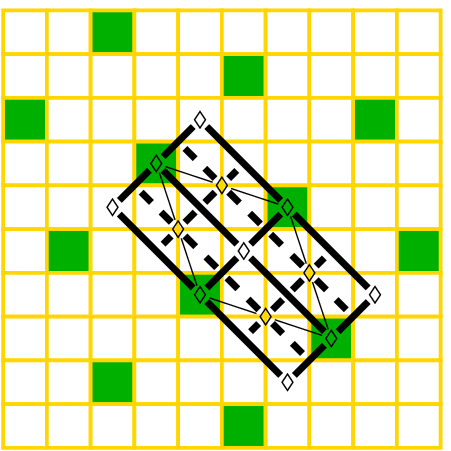}\hskip 5 pt \epsffile{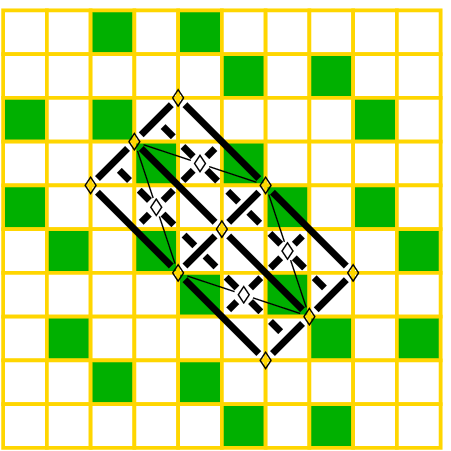}

\noindent (a)\hskip 1.66 in (b)
\smallskip


\noindent Figure 13. Fabrics illustrating two species with geometrically isomorphic%
\donote{2}{This is the stronger geometrical sense of isomorphism used in [Thomas, MS, p.~6] and explained in [Gr\" unbaum and Shephard, 1987, p. 38].}
 groups of Roth type 28. 

\noindent a.~The (8, 3) satin (species $28_e)$. \hskip 10 pt b.~8-5-1 (species $28_n$).
\smallskip
\smallskip\noindent\epsffile{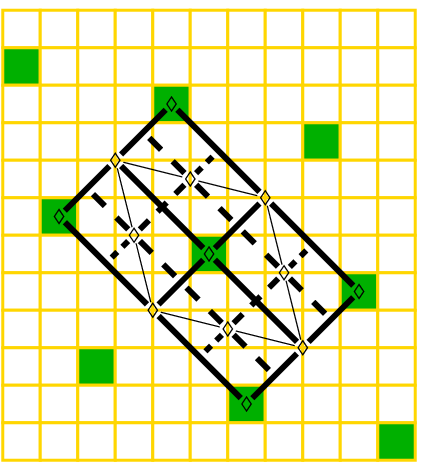}\hskip 10pt \epsffile{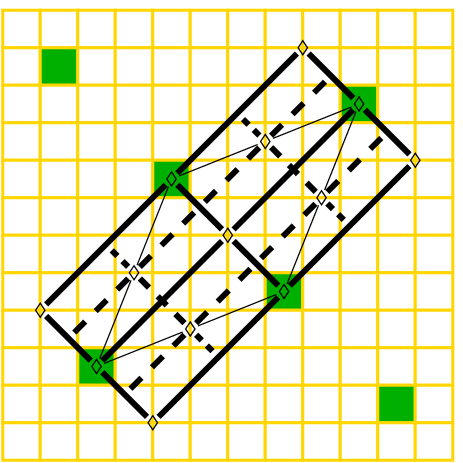}

\noindent (a)\hskip 1.65 in (b)
\smallskip


\noindent Figure 14. Species-$28_o$ examples. \hskip 10 pt a.~The (15, 11) satin with central rectangle $3\beta$ by $5\beta$. \hskip 10 pt b.~The (21, 13) satin with central rectangle $3\beta$ by $7\beta$.
\smallskip

If the central rectangle has corners at cell centres, and standard isonemal spacing to distribute them, odd-odd dimensions in $\delta$ units make all the nine-configuration be at cell centres too and there are no strand-to-strand half-turns.
This must be ruled out.
Odd-even dimensions again give the appropriate mixture of sorts of centres, essentially reversing the roles played in species $28_e$.
Figure 13b illustrates an example 8-5-1 of this kind with $a=1$, $b=2$, period 8.
We shall call this species $28_n$ 
because it is not exemplified by satins; its period and order are $4ab$ divisible by 8.
Roth mentions [1995, p.~322] that 8-1-1 and 8-5-1 differ in their groups in spite of their groups being of the same type.
Yet a third sort of location is possible.
Supplementary half-turn centres located in cell boundaries must appear matched diagonally across the central rectangle because they are related by a half-turn, whereas adjacent corners are related by reflection and so must be opposite in the side-top/bottom dichotomy (Figure 14).
The sides of such rectangles are odd in $\beta$ units.
The centre can be either at a cell centre (Figure 14a) or cell corner (Figure 14b).
Because the quarters of the nine-configuration have the same dimensions as the central rectangle, the configuration's corners will be in the same position (cell corner or cell centre) as the centre of the central rectangle, and the mid-sides (the corners of the rhombic $G_1$ lattice unit) will be in the other position.
All that is required to distribute the half-turn centres is that the central rectangle have standard isonemal spacing in $\beta$ units, making the whole nine-configuration have odd-odd isonemal spacing in $\delta$ units ($a=1$ included, as the 2/1 twill with central rectangle $\beta$ by $3\beta$ illustrates).
More interesting examples of both sorts are the $(15, 11)$ satin with lattice-unit centre at a cell centre and the $(21, 13)$ satin with lattice-unit centre at a cell corner, both illustrated in Figure 14.
(For uniformity, all satin offsets are given as odd.)
But note that the rhombs in Figure 14 could be drawn differently, centred on the dark cell in Figure 14b and with corners at the dark cells in Figure 14a.
The third arrangement is not subdivided in an essential way into two but is only a single species $28_o$ (for odd). 
Orders of this species are $ab$, the only odd orders in the prefabrics with perpendicular axes.

Four of the five spacings just determined, $27_e$, $27_o$, $28_e$, and $28_o$, will recur in Roth types 29 to 32 respectively.
No analogue of $28_n$ can exist---another reason for the n---because it would place \diab in cells.
\smallskip
\noindent{\bf 8. Rhombic satins}

\noindent Rhombic satins are the simplest examples of Roth type 28.
Two of its subtypes have been illustrated by them, $28_e$ by even-order rhombic satins and $28_o$ by odd-order rhombic satins.
Those of even order require that either the cells containing the corners of the rhombic lattice unit be darkened or the centre since both fall in cells as in the $(8, 3)$ satin in Figure 13a.
The nine-configuration could be displaced by half of its length or width.
Those of odd order likewise can have the rhombs drawn in either way---again by displacement of the nine-configuration.
The even and the odd satins, however, are essentially distinct.
The family of even rhombic satins like $(8, 3)$ begin by being based on a central rectangle with corners havingtwill spacing and then with rectangular-satin spacing with central-rectangle corners at cell corners.
The 3/1 twill gives the $(8, 3)$ satin, the 7/1 twill the $(16, 7)$ satin, the 11/1 twill the $(24, 11)$ satin, and so on.
Then the $(12, 5)$ satin gives the $(24, 5)$ satin, the $(20, 9)$ satin gives the $(40, 11)$ satin, and so on.
The family of odd rhombic satins like $(15, 11)$ does not begin with central rectangle corners at twill spacing because, as the 2/1 twill illustrates, that spacing with corners on cell boundaries gives rise to twills not satins.
The smallest central rectangle---that giving the $(15, 11)$ satin---is $3\beta$ by $5\beta$, which would be the spacing of the $(30, 11)$ satin if it were in $\delta$ units.
Likewise the next such rhombic satin, $(21, 13)$ has central rectangle $3\beta$ by $7\beta$, which would be the spacing of the $(42, 13)$ satin if it were in $\delta$ units.
\smallskip
\settabs \+ \hskip 0.5 in & \hskip 1.5 in & \hskip 0.9 in &\cr
\+& Rectangular        &Rhombic   &Rectangular\cr
\+& (12, 5)            &(24, 5)\cr
\+& (20, 9) = (20, 11) &(40, 11)  &(40, 9)\cr
\+& (24, 7)            &(48, 7)\cr
\+& (28, 13)           &(56, 13)\cr
\+& (30, 11) = (30, 19)&          &(60, 11), (60, 19)\cr
\+& (36, 17) = (36, 19)&(72, 19)  &(72, 17)\cr
\+& (40, 9)            &(80, 9)\cr
\+& (42, 13) = (42, 29)&          &(84, 13), (84, 29)\cr
\+& (44, 21)           &(88, 21)\cr
\+& (48, 17)           &(96, 17)\cr
\smallskip
\noindent Table 1. Pairing (in columns 1 and 2) of rectangular and rhombic satins. Column 3 lists rectangular satins that might be expected to be rhombic.
\smallskip
This analysis sheds some light on what little is known in the theory of satins.
We know [Gr\"unbaum and Shephard, 1980, Theorem 3] that an $(n, s)$ satin is isonemal if and only if $s^2 \equiv \pm 1$ (mod $n$).
(The condition $s^2 \equiv -1$ (mod $n$) is the condition that the $(n, s)$ satin be square.
Square satins do not have perpendicular symmetry axes.)
And [Gr\"unbaum and Shephard, 1980, Theorem 4], as cited in Section 4, a non-square isonemal satin is rectangular if $n$ is even and $s^2 \equiv 1$ (mod $2n$). 
Otherwise it is rhombic.
So it is clear that, corresponding to each rectangular satin $(2m, s)$ with $s^2 \equiv 1$ (mod $4m$), there is a rhombic satin $(4m, s)$.
This condition pairs rectangular and rhombic satins as in Table 1 (columns 1 and 2).
Several satins (column 3), which would be accounted for if rhombic, are rectangular.
Numerically, the reason is that the offset $s$ of the $(n, s)$ satin is large enough that $s^2 \equiv 1$ (mod $4n$), so that the $(2n, s)$ satin meets the rectangular-satin condition.
In several cases, $s^2 -1 = 4n$ exactly.
The geometrical explanation is that, while the $(30, 11)$ and $(42, 13)$ satins are rectangular, their rectangles, 30 and 42 not being divisible by 4, are odd by odd (Theorem 1), and they are not used to produce rhombic satins of even order but rather of odd order.
The entries that belong in the gaps in column two of Table 1 are $(15, 4) = (15, 11)$ and $(21, 8) = (21, 13)$, so much smaller because their odd-by-odd central rectangle is measured in $\beta$ units.

There is no reason to place much emphasis on the relation between rectangular and rhombic satins since either can easily be turned into one of the other sort and of larger period by deletions (of alternate mirrors in the rhombic case) and the process is reversible if certain divisibility conditions are fulfilled.
\smallskip
\noindent {\smc Theorem 2.} Rhombic satins of order odd and even fall in species $28_o$ and $28_e$ respectively.
\medskip
\noindent {\bf 9. Crystallographic type} {\bfsl cmm} {\bf with side-preserving glide-reflec\-tions: Roth types 29--32}

\noindent {\it 29.} The alternative chosen at the beginning of Section 7 was to have side-reversing glide-reflections between the mirrors.
The other alternative is side-preserving glide-reflections between the mirrors.
Then the half-turns with centres at their intersections must be side-reversing as in species 19.
The half-turns (and mirrors) are absent from the side-preserving subgroup, but the glide-reflections are there, making it be of crystallographic type $pgg$ with the rectangle of the nine-configuration as the lattice unit.
While much of the above discussion of type 27 carries over, care must be taken that \diaas at mirror intersections and \diabs at axis intersections consistently appear either together on every strand boundary or separately on alternate boundaries to keep their effects on the strands consistent, i.e., for isonemality.
If \dia and \diab appear only at cell corners and never on the same side of a strand, there are no implications for strand symmetry, the subgroups are of type $11/-$, and the Roth type is 29.
Because of where they are constrained to intersect, the axes are in mirror position.
The smallest dimensions, $2a\delta$ by $2b\delta$, of the central rectangle, which is the same size as one quarter of the $H_1$ ($pgg$) lattice unit, are $2\delta$ by $4\delta$ as with the even-even species $27_e$.
Roth's example (Figure 15a) is smallest possible with order 16 equal to one quarter of the\break
\smallskip
\noindent\epsffile{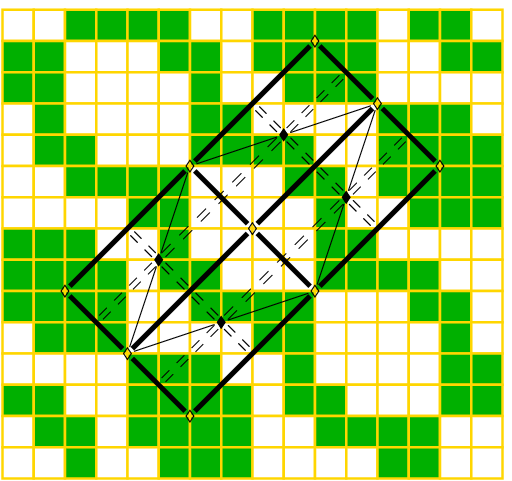}\hskip 10pt \epsffile{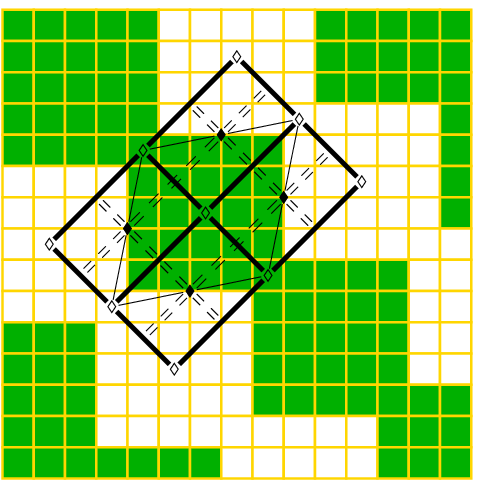}

\noindent (a)\hskip 2 in (b)
\smallskip
\noindent Figure 15. a. Roth's example 16-2499 of species 29. \hskip 10 pt 

\noindent b. Example of species 31 of order 24.

\vfill\eject

\noindent period (the size of the $H_1$ lattice unit) since each portion of strand of order length appears in four guises in successively adjacent strands on account of being transformed alternately by \dia and \diabb .
The genus is pure V, which prevents $H_1$ from being transitive on strands, and the order is $8ab$, always divisible by 16.

\noindent {\it 30.} To produce strand subgroups of type $11/11$ and so Roth type 30, we need the \dia and the \diab lattices to fall at cell corners but both to fall on every strand boundary.
The lattice of \diaas and the lattice of \diabs (both with the same spacing) must be distributed to every strand boundary, and the spacing that will do that has already been considered in the discussion of species $27_o$.
Because $H_1$ has \dia on every strand boundary, it is transitive on strands.
Roth's examples, 12-189-2 and 12-315-4, with $a=1$ and $b=3$ are smallest possible and confirm that $a=1$ is allowed.
An example not using $a=1$ as a dimension has $a=3$, for which the smallest possible $b$ is 5, and is given in Figure 16. %
The period is 120 and order 60.
Order in general is $4ab$, divisible by 4.

\smallskip\noindent\epsffile{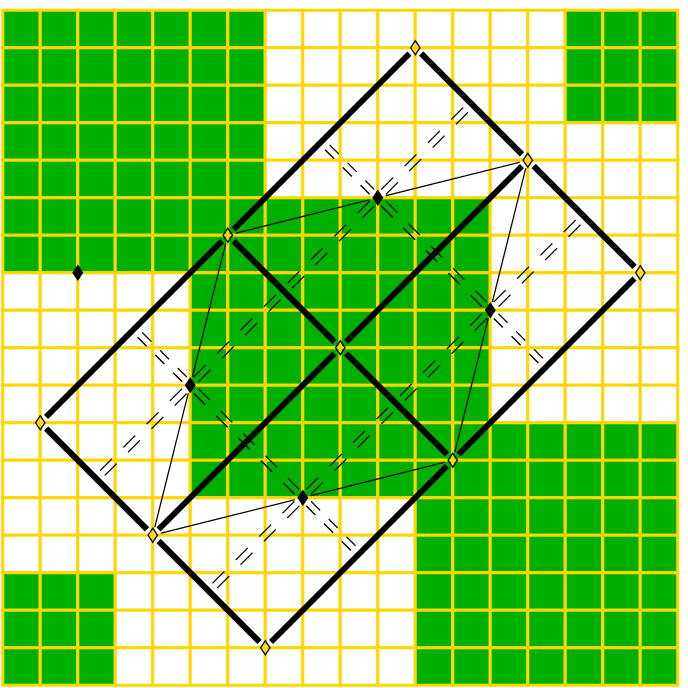}
\smallskip


\noindent Figure 16. Example of species 30 of order 60. The single \diab illustrated outside the $H_1$ lattice unit is on the same edge of a strand as a \dia in a corner of the lattice unit.
\smallskip

\noindent {\it 31.} To produce strand subgroups of type $12/-$ and so Roth type 31, we need the \diaas to fall at cell centres and the supplementary \diabs at cell corners.
Odd-even spacing as in species $28_e$ for the central rectangle is what is required ($a=1$ allowed) so that its centre falls in a cell.
Roth's examples both have $a=1$, $b=2$ but this is unnecessary, the next allowable sizes having $a=1$, $b=4$, and $a=3$, $b=2$, the latter illustrated in Figure 15b. 
Period is 48 and order 24; the period is larger than the order because the \diab but not four times because of the strand-symmetry-producing {\sym\char'006}s are in strands not between them.
Order in general is $4ab$, divisible by 8.
As in types 22 and 23, the even-by-even dimensions of the lattice unit of $H_1$ prevent it from being transitive on strands without \diaa , but they are in cells and so cannot help.
$H_1$ of type $pgg$ is not transitive on strands.

\noindent {\it 32.} Last of all, to produce strand subgroups of type $12/12$ and so Roth type 32, we need both \dia and \diab to fall in strands.
This requires the spacing of species $28_o$ with $a=1$ allowed (as in Roth's examples, both twills) in $\beta$ units.
It must be the \diabs that fall on the cell boundaries because only the \diaas can fall at cell centres.
The order in general is $2ab$ rather than just $ab$ as for species $28_o$ on account of the {\sym\char'007}s, which in strands make order even.
The smallest central rectangle not having a dimension 1 is $3\beta$ by $5\beta$, the spacing, cast into $\beta$ units, of the $(30, 11)$ satin.
The four species 29--32 are all based on twills or rectangular satins in this way.
With these dimensions, the designs available have six cells specifiable arbitrarily, but the $2^6$ designs, being self-complementary on account of the \diabb , all appear twice.
The $2^5$ distinct designs based on this group, which is a subgroup of several groups with additional symmetry, include designs with larger symmetry groups: the 5/5 twill, the 1/3/1/3/1/1 twill, and the 2/2/1/2/2/1 twill with a $\beta$ by $5\beta$ central rectangle and order 10, the 3/3 twill with a $\beta$ by $3\beta$ central rectangle and order 6, and plain weave of order two, which can be seen as a design like this with a $\beta$ by $\beta$ central rectangle.
The other $2^5 - 5 = 27$ designs are pleasing tessellations, sometimes surprising.
Half a dozen examples are given in Figure 17.
\smallskip
\noindent
\epsffile{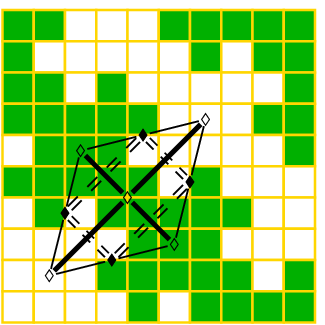}\hskip 10 pt\epsffile{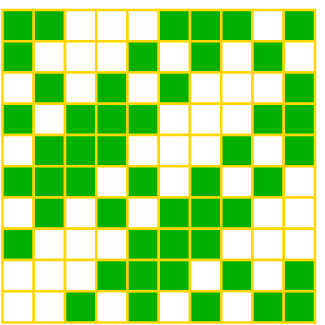}\hskip 10 pt \epsffile{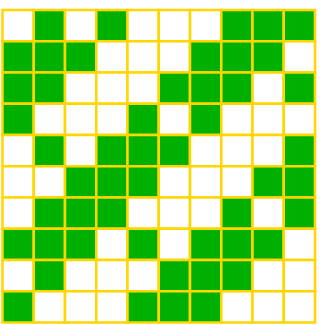}

\noindent
\epsffile{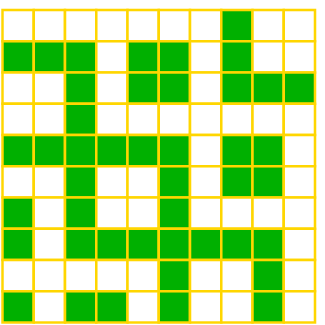}\hskip 10 pt\epsffile{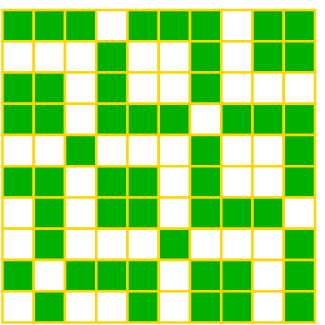}\hskip 10 pt \epsffile{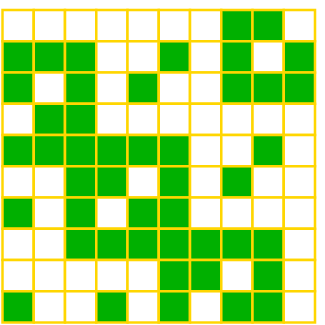}
\smallskip
\noindent Figure 17. Six designs of species 32 and order 30.

\smallskip

The parameters $a$ and $b$, chosen to ensure isonemality and determining period and order, are throughout dimensions in $\delta$ or $\beta$ units of some rectangle that needs to be distributed over the plane.
When it is clear what the basic rectangle is and its corners are at cell corners, the dimensions are in $\delta$ units (species---in groups with similar spacing---11, 13, and 15; 17, 19, and 21; 25, $27_e$ and 29; $27_o$ and 30; $28_n$). 
Likewise, when mid-lines of strands are involved, dimensions are in $\beta$ units (grouped species 12, 14, and 16; 18 and 20; 26; $28_e$ and 31; $28_o$ and 32).
There are only three Roth types that do not fit into these patterns, 22, 23, and 24, each of which has an $a\delta$-by-$b\delta$ rectangle distributed; it is just less obvious what rectangle that is.
\smallskip
\noindent {\bf 10. Doubling} 

\noindent The question was raised in [Thomas, MS] which designs can be doubled and remain isonemal as plain weave is doubled to 4-3-1, each strand being replaced by a pair of strands with the same behaviour.
Doubling was introduced into the weaving literature by Gr\"unbaum and Shephard in their [1980] applied to square satins like 5-1-1.
Since doubling can no more be expected to preserve isonemality twice than it can be expected that three strands can be set to the work of one, the result of doubling, when it can be done preserving isonemality, will always be designs that cannot be doubled.
We have two theorems, both of which apply to all species 1--32, the latter applying trivially to species 1--10.
\smallskip
\noindent{\smc Theorem 3.} No non-exceptional isonemal prefabric with symmetry axes and strand subgroup of type $11/-$ or $11/11$ can be doubled and remain isonemal.
\smallskip
\noindent{\it Proof.} The mapping between the two strands that replace each strand needs to transform the one strand to the other by translation, reflection, or a half-turn.
No translation is possible, for if, say, the lower half of a strand is translated to the half above, then the half above is moved by the same translation to the lower half of the strand above. 
The latter movement (along the strand length) can only be by the offset between the strands and be a symmetry of the design, and that number is relatively prime to the order, which is the minimum movement (again along the strand length) mapping the halves of the strand one to the other.
No suitable reflections are available in the Roth types 1 to 32. 
And the presence of the needed half-turns makes the strand subgroup be of type $12/11$, $12/12$, or $12/-$.
No such mapping is possible.
\medskip
\noindent{\smc Theorem 4.} All non-exceptional isonemal prefabrics with symmetry axes and strand subgroup of type $12/11$, $12/12$, or $12/-$ can be doubled and remain isonemal.
\smallskip
\noindent{\it Proof.} The species to which this theorem applies are 12, 14, 16, 18, 20, 23, 24, 26, 28, 31, and 32, and when doubled they become prefabrics of type respectively 11, 13, 15, 17, 19, 21, 22, 25, 27, 29, and 30 because centres of half-turns on mid-strand lines become half-turn centres at cell corners in the doubled prefabric mapping each pair of replacement strands, both warps and wefts, to each other. Other strand-to-strand mappings remain.
\smallskip
The images of all of these transformations under doubling are in species that cannot be doubled, each of which contains a core of designs that are images under doubling. 
\smallskip
\noindent {\bf 11. Halving} 

\noindent One would expect that what had been doubled could be halved although the reverse might not be true since the operations are not really inverses.
Halving is first the removal of every other strand in each direction, making an intermediate construction called a pseudofabric by Gr\"unbaum and Shephard [1988], and then the widening of all the strands uniformly to produce another prefabric.
So halving undoes doubling, but doubling does not normally reverse halving.
In fact, we have the extension to these prefabrics of the theorem on this topic in [Thomas, MS].
\smallskip
\noindent {\smc Theorem 5} Isonemal prefabrics of species 11 to 32, except those satisfying the two conditions (1) having only glide-reflections with axes not in mirror position in both directions (species $12_o$ and $14_o$) and (2) having both lattice-unit dimensions greater than $2\delta$, can be halved with only isonemal results.
\smallskip
For the proof we require the cells of the plane to be numbered in fours: group the cells into square blocks of four that tessellate the plane as do the cells, and number the cells of each block as the quadrants of the Cartesian plane are numbered.
Beginning with a prefabric, one can then produce any half-fabric or {\it factor} by preserving the crossings at cells assigned any one of the four numbers and discarding the rest.
To prove the theorem we need the following lemma from [Thomas, MS].
\smallskip
\noindent {\smc Lemma 2.} The action of a reflection or glide-reflection with axis in mirror position on the numbering of the cells is displayed in Table 2, where the {\it parity} of an axis is the parity of the cells through which it runs, {\it length} refers to dimension of the $G_1$ lattice unit parallel to a mirror in $\delta$ units, and {\it interchange of numbers} means specifically the interchange of 1 and 3 or 2 and 4.

\medskip

\noindent {\it Proof of Theorem 5.} 
If a lattice-unit dimension is $2\delta$ or less, then the prefabric is a twill and can be halved to another twill, all of which are isonemal.
It remains to show that isonemal prefabrics not in species $21_o$ or $14_o$ and having both lattice-unit dimensions greater than $2\delta$ can be halved with idsonemal results.

We shall work though the types considering what does or does not allow halving.
Half-turns at cell corners between all adjacent strands, since such half-turns preserve each parity and interchange numbers with each parity, act on alternate strands to preserve numbering in each parity; a\break
\smallskip

\noindent\hskip 0.75 in Parities 

\noindent\hskip 0.75 in of axis \hskip 0.5 in Same \hskip 1.2 in Opposite

\noindent\hskip 0.75 in and cells

\vskip 2 mm

\noindent Operation
\vskip 2 mm

\settabs \+ \hskip 1.25 in & \hskip 1.6 in & \cr

\+Odd glide & interchanges numbers & preserves numbering\cr

\smallskip

\+Even glide & preserves numbering &interchanges numbers\cr

\smallskip
\+Reflection & preserves numbering & interchanges numbers\cr

\+(even length)\cr

\smallskip
\+Reflection & preserves numbering & preserves numbering when\cr

\+(odd length) && combined with translation\cr

\smallskip
\noindent Table 2. Effects of reflections and glide-reflections on cell numbers.
\medskip

\noindent cell-corner half-turn interchanges 2 and 4, 1 and 3, but two successive such half-turns preserve 1s, 2s, 3s, and 4s.
These half-turns are present in species 11, 13, 15, 17, $18_s$, 19, 21, 22, 23, 25--32.
All that is needed to make such halved designs isonemal is a warp-to-weft symmetry transformation.
These are shown to exist in the next six paragraphs; the remaining species ($12_e, 14_e, 16, 18_o, 18_e, 20$, and 24) are then discussed in a further three paragraphs.
Types 11, 13, and 15 can be treated together. 
They have mirror-position axes of glide-reflection odd by odd distances apart in $\delta$ or odd by even. 
These meet at cell corners, respectively cell centres, meaning that they lie in cells of both parities, respectively just one parity $p$. 
In the former case, the perpendicular glides are both odd and so preserve the numbering of cells of both parities by Lemma 2. 
In the latter case, the perpendicular glides are of parities odd and even and so, lying in cells of the same parity $p$, transform warps to wefts while preserving the numbering of both parities of cell $\overline p$ and $p$ respectively.

Types 17, 19, 21, and 22 can be dealt with together by discussing 17 because species $19_o$ and 22 behave like $17_o$, and $19_e$ and 21 like $17_e$.
They have mirrors crossing axes of glide-reflection half-way between the centres of half-turns along the axes.
When the distance between the centres along the axes is odd in $\delta$, that is, in species $17_o$ and on one orientation of the basic rectangle in $17_e$, the mirror and axes intersect in a cell centre, and the mirror and axes run through cells of the same parity $p$.
But in that case the glide is odd, and so the mirrors preserve numbering in cells of parity $p$ and the glide-reflections preserve numbering in cells of parity $\overline p$, giving the warp-to-weft transformations required.
When the distance between the centres along the axes is even, that is in the other orientation of basic rectangles in species $17_e$, the mirror and axes intersect at a cell corner, and the mirror and axes run through cells of opposite parities.
But in that case the glide is even, and so the mirrors and the glide-reflections preserve the numbering in cells of opposite parities, giving the warp-to-weft transformations required.
These considerations of the two ways of orienting these basic rectangles was adumbrated in the second paragraph on species 15.

Species $18_s$, 23, $28_o$, and 32 have axes of glide-reflection an odd multiple of $\beta$ apart therefore lying in cells of opposite parity. 
They serve alternately to transform factors of both parities.
Because the dimension of the lattice unit of $18_s$ parallel to the mirrors is odd in $\delta$, the mirrors too transform factors of both parities as in species 6 in [Thomas, MS].
Because they are an odd multiple of $\beta$ apart, the mirrors of type $28_o$ in {\it both} directions alternately transform factors of both parities.

Species 25, 27, $28_n$, 29, and 30 have mirrors in both directions a multiple of $\delta$ apart and so each direction in only one parity, but the two directions are in cells of different parities, for they meet at cell corners, and so preserve numbering in the cells of each parity.

Species $26_o$ and $26_e$ have mirrors at half an odd multiple of $\delta$ apart in one direction ($26_e$) or both ($26_o$) and so in cells of both parities, therefore preserving numbering in both parities.

Species $28_e$ and 31 have central rectangles with corners at cell corners and with sides odd and even in $\delta$.
Mirrors through their central cells bisect their sides, the even sides being bisected at cell corners.
An axis of even-glide glide-reflection and a mirror, meeting at a cell corner, lie in cells of opposite parity; each preserves the numbering of cells of its own parity, that is of cells of opposite parities, giving the warp-to-weft transformations required, as in the latter case of $17_e$.

Then there are the species that do not have centres of half-turns at cell corners along all strand boundaries.
Species $12_e$, $14_e$, and 16 have axes of glide-reflection odd $a$ by even $b$ apart in $\beta$, with only those $a\beta$ apart in mirror position and with $(a, b)=1$.
The spacing of the axes in mirror position ensures that they fall in cells of both parities alternately, and so one or the other set preserves the numbering of cells of each parity.
Each factor of each parity has alternate glide-reflection axes in it.
The lattice unit of each factor has width the distance from one axis to that four axes away, $4a\beta$, cut in half by the halving, and so $2a\beta=a\delta$.
Its length is that of the lattice unit $2b\beta$, cut in half by the halving to $b\beta=b\delta/2$.
The result is an isonemal design of species $2_m$ with axes of side-reversing glide-reflection in mirror position and with lattice-unit dimensions the relatively prime $a$ and $b/2$ in $\delta$ units as required.
The symmetry group could be larger and contain one of Roth type $2_m$ properly  since additional symmetry can easily be created by the deletion of cells.
(Species $12_o$ and $14_o$, however, have axes odd by odd distances apart in $\beta$ with none in mirror position, consequently of no use for transforming warps to wefts while preserving parities, still less numbers.)

The glide-reflection axes of species $18_o$ and $18_e$ are not in mirror position and so are of no help.
The mirrors perpendicular to them in both species are spaced an odd multiple $a\beta$ apart.
The same is true of species $20_o$ and $20_e$, $24_o$ and $24_e$ (there is no analogue of $18_s$).
Much the same thing happens for species $18_e$, $20_e$, and $24_e$ as happens in species $12_e$, $14_e$, and 16 but with mirrors instead of axes.
I shall discuss $18_e$.
The mirrors fall alternately in cells of each parity.
Each factor of each parity has alternate mirrors in it.
The lattice unit of each factor has width the distance from one mirror to that four mirrors away, $4a\beta$, cut in half by the halving, and so $a\delta$.
Its length is the length of an $18_e$ lattice unit, $2b\beta$, cut in half by the halving to $b\delta/2$.
The result is an isonemal design of species 5 with lattice-unit dimensions the relatively prime $a$ and $b/2$ in $\delta$ units as required, regardless of the parity of $b/2$.
Again the symmetry group might only {\it contain} one of Roth type 5.

Somewhat the same thing happens for species $18_o$, $20_o$, and $24_o$, again with mirrors instead of axes of glide-reflection.
I shall discuss $18_o$.
The lattice unit of each factor has width the distance from one mirror to that four mirrors away, $4a\beta$, before being cut in half by the halving to $a\delta$.
Its length is the length of two $18_o$ lattice units since it has to be even in $\delta$, $4b\beta$, cut in half by the halving to $b\delta$.
The mirrors fall alternately in cells of each parity, but both serve to transform factors of both parities.
Each factor of each parity has alternate mirrors in it and alternate mirrors missing from it.
Because of the odd length of the species-$18_o$ $G_1$ lattice unit, the missing mirrors act as axes of side-reversing glide-reflection not in mirror position by Lemma 2 (bottom right entry of the table).
The result is an isonemal design of species $8_o$ with lattice-unit dimensions the odd relatively prime $a$ and $b$ in $\delta$ units as required.
$20_o$ is much the same, but because of the $\tau$-translation between adjacent $G_1$ lattice units the result is of species 10 rather than $8_o$.
Again it is not quite accurate to say that $G_1$ is of type $8_o$ or 10, because while the design has those symmetries the actual group may be a larger group that contains the group of the type predicted as a subgroup.
Since all that is needed is a transitive group of symmetry operations, in all cases the theorem is proved.
\smallskip
The half-turns that are a feature of these designs have effects on the halved patterns beyond that mentioned at the beginning of the above proof.
If the centre of a half-turn without $\tau$ is at a cell corner, then the two factors of each parity are isometric images of each other (the same).
The species listed at the beginning of the proof included designs with \diab at corners.
Those with \dia at a corner are 11, 13, 17, $18_s$, 21, 22, 25, 26, 27, $28_e$, $28_n$, 29, 30, 32.
If the \diab is at a cell corner, then the two factors of each parity are isometric images of each other with side reversal (complementary).
Those with \diab at a cell corner are 15, 19, 21, 22, 23, 29, 30, 31.
If a \dia is on a cell edge, then the factors of each parity together are isometric images of the factors of the other parity (the same in pairs).
Those with \dia on a cell edge are 12, 14, $18_o$, $18_e$, 24, $28_o$.
If a \diab is on a cell edge, then the factors of each parity together are isometric images of the factors of the other parity with side reversal (complements in pairs).
Those with \diab on a cell edge are 16, 20, 24, 32.
As the lists above indicate, these are not exclusive conditions in species numbered beyond 20.
\smallskip
\noindent {\bf 12. Fabrics of a given order} 

\noindent As in [Thomas, MS], the information dispersed in sections 3, 5--7, and 9 allows all fabrics of a given order to be found because the symmetry groups for a given order can be determined. 
It is not always useful to divide prefabrics as finely as into species.
In \S 10 it was not.
One place where such division is useful is in this determination.
If the order is odd, then {\it among these species} (a restriction that persists throughout the section), it must be $ab$, the order of a fabric of species $28_o$ with parameters $a$ and $b$, whether $ab$ is composite or prime, since $a=1$ is allowed.
If the order is even, then the possible parameters and species are displayed in Table 3.
Various factorizations need to be used, and the allowance or disallowance of $a=1$ and of using lattice units both ways (these two ways could determine subsubspecies) need to be taken into account.
It is now even easier than at the end of [Thomas, MS] to see why fabrics of order 16 have not been catalogued.
\smallskip
\settabs\+ Condition\hskip 10 pt & Order\hskip 10 pt &  Species\cr
\+ Condition &Order & Species\cr
\+ $2 | N$, $4\!\not |\, N$ & $2ab$  & $12_o$, $14_o$, $18_o$, $20_o$, $26_o$, $27_o$, 32\cr
\smallskip
\+ $4 | N$                 & $2ab$  & $12_e$, $14_e$, $16$, $18_e$, $18_s$, $20_e$, $26_e$\cr
\smallskip
\+                          & $4ab$  & $11_o$, $13_o$, $15_o$, $17_o$, $19_o$, $23_o$, $24_o$, $25_o$, 30\cr
\smallskip
\+ $8 | N$                  & $4ab$  & $11_e$, $13_e$, $15_e$, $17_e$, $19_e$, 21, 22, $23_e$, $24_e$,\cr
\+&& $25_e$, $28_e$, $28_n$, 31\cr
\smallskip
\+ $16 | N$                 & $8ab$  & $27_e$, 29\cr
\smallskip
\noindent Table 3. Species in which prefabrics of even orders $N$ may be found. 
\smallskip

From the symmetry groups as determined in the previous paragraph, the designs can be produced by the usual method [Thomas, MS, Section 5], which produces---as well as what is wanted in any particular case---also designs with symmetry groups having the desired symmetry group as a subgroup as in [Thomas, MS] and also prefabrics that fall apart.
A difference between species 1--10 and species 11--32 is that, as well as falling apart as a special case of halving, another possibility emerges (next paragraph).
To deal first with the analogue of what happens in [Thomas, MS], Roth [1993] pointed out that prefabrics of species 15, 19, 23, and 31 had the genus (IV in the cases of 15 and 19) or genera (II and IV in the other cases) to fall apart with alternate warps and wefts lifting off the others.
Where the warps and wefts that will lift off cross over the others, half of the cells of the design are predetermined as dark or pale. 
The restrictions on where through this array axes and mirrors can run are those mentioned in [Thomas, MS]:
Along diagonals of predetermined colouring can be placed mirrors, axes of side-preserving glide-reflection with odd glides, and axes of side-reversing glide-reflection with even glides.
Between those diagonals, that is, through cells not predetermined, can be placed axes of side-preserving glide-reflection with even glides and axes of side-reversing glide-reflection with odd glides.
The orders of designs of species $15_o$, $19_o$, and $23_o$ are always divisible by 4, those of species $15_e$, $19_e$, $23_e$, and 31 by 8.
The smallest actual example of $15_e$, as mentioned in section 5, is precisely one of these prefabrics that falls apart, 8-5-2*, and they continue to occur, e.g. 16-85-2* and 16-277-2*, which use the same rectangle both ways.
Species $15_o$ designs begin with 12-21-4* and 12-69-4* using the same rectangle both ways.
The smallest examples of $19_o$ are 12-21-3* and 12-69-3*, which use the same basic rectangle ($a=1$, $b=3$) differently oriented, and of $19_e$ 16-85-3* and 16-277-3* with the same basic rectangle ($a=1$, $b=4$) differently oriented.
The smallest actual example of $23_o$ is 4-1-1*, continuing with 12-21-1* and 12-21-2*.
The smallest example of $23_e$ is 8-5-1*, continuing with 16-85-1* and 16-1105-1*.
Examples of species 31 begin with 8-5-3* and continue with 16-85-5* and 16-1105-3*.

As well as falling apart that is a species of halving as they have both been described above, isonemal prefabrics can also fall apart in a different way with adjacent pairs of warps and of wefts lifting off alternate adjacent pairs.
This phenomenon has been discussed in [Hoskins and Thomas, 1991] and [Roth, 1995].
It can occur in only two species among those considered here, 21 and 29, of pure genus V with orders divisible by 8 and 16 respectively.
The rules for falling apart are similar (and imply only that orders be divisible by 4) except that they apply to adjacent pairs of strands:
Warps corresponding to predominantly (3/4) dark columns together with wefts corresponding to predominantly (3/4) pale rows can be lifted off the others, i.e., warps corresponding to predominantly (3/4) pale columns and wefts corresponding to predominantly (3/4) pale rows.
Half of the cells of such a design are where warps and wefts that can be lifted off cross over what is left behind; accordingly half of the cells are predetermined as dark or pale in square blocks of four cells.
Along the diagonal of a block of predetermined colouring can be placed a mirror, an axis of side-preserving glide-reflection with glide twice an odd number, or an axis of side-reversing glide-reflection with glide twice an even number.
Odd glides are impossible, as are mirrors and axes not along diagonals of blocks of predetermined colouring or blocks not predetermined.
Along the diagonals of blocks of cells not predetermined in colour, that is between blocks of cells of predetermined colour, can be placed an axis of side-preserving glide-reflection with glide twice an even number or an axis of side-reversing glide-reflection with glide twice an odd number.
Species 21 examples begin with 8-9-1* and 8-3-1*, which is 4-1-1* doubled, and continue with 16-153-1* and 16-51-1*, which is 8-5-1* doubled.
Species 29 examples begin with 16-51-2*, which is 8-5-3* doubled, and 16-153-2*.
Plainly they must be rooted out of any list that purports to consist only of designs of fabrics.

The results of sections 10--12 need to be extended to the remaining prefabrics.
\medskip
\noindent {\bf References}
\smallskip

\hangindent 12 pt\noindent Branko Gr\"unbaum and Geoffrey C.~Shephard, Satins and twills: An introduction to the geometry of fabrics, {\it Mathematics Magazine} {\bf 53} (1980), 139--161.

\hangindent 12 pt\noindent ------ A catalogue of isonemal fabrics, in {\it Discrete Geometry and Convexity}, Jacob E.~Goodman {\it et al.,} eds. {\it Annals of the New York Academy of Sciences} {\bf 440} (1985), 279--298.

\hangindent 12 pt\noindent ------ An extension to the catalogue of isonemal fabrics, {\it Discrete Mathema\-tics} {\bf 60} (1986), 155--192.

\hangindent 12 pt\noindent ------ {\it Tilings and patterns,} W.H.~Freeman, New York, 1987.

\hangindent 12 pt\noindent ------ Isonemal fabrics, {\it American Mathematical Monthly} {\bf 95} (1988), 5--30.


\hangindent 12 pt\noindent J.A.~Hoskins and R.S.D.~Thomas, The patterns of the isonemal two-colour two-way two-fold fabrics, {\it Bull.~Austral.~Math.~Soc.}~{\bf 44} (1991), 33--43.

\hangindent 12 pt\noindent Richard L.~Roth, The symmetry groups of periodic isonemal fabrics, {\it Geometriae Dedicata} {\bf 48} (1993), 191--210.

\hangindent 12 pt\noindent ------ Perfect colorings of isonemal fabrics using two colours, {\it Geometriae Dedicata} {\bf 56} (1995), 307--326.

\hangindent 12 pt\noindent D.~Schattschneider, The plane symmetry groups: Their recognition and notation, {\it American Mathematical Monthly} {\bf 85} (1978), 439--450.

\hangindent 12 pt\noindent R.S.D.~Thomas, Isonemal prefabrics with only parallel axes of symmetry, {\it Discrete Mathematics}, to appear. At www.arxiv.org/math.CO/0612808.

\bye